\documentclass[reqno,a4paper]{amsart}
\usepackage{amssymb}
\usepackage{amsmath}
\begin{document} 
\newtheorem{prop}{Proposition}[section]
\newtheorem{Def}{Definition}[section] \newtheorem{theorem}{Theorem}[section]
\newtheorem{lemma}{Lemma}[section] \newtheorem{Cor}{Corollary}[section]

\title[Gross-Pitaevskii equation]{\bf Unconditional global well-posedness for 
the
3D Gross-Pitaevskii equation for data without finite energy}
\author[Hartmut Pecher]{
{\bf Hartmut Pecher}\\
Fachbereich Mathematik und Naturwissenschaften\\
Bergische Universit\"at Wuppertal\\
Gau{\ss}str.  20\\
42097 Wuppertal\\
Germany}
\email{\tt pecher@math.uni-wuppertal.de}
\date{}

\begin{abstract}
The Cauchy problem for the Gross-Pitaevskii equation in three space dimensions 
is shown to have an unconditionally unique global solution for data of the 
form $1 
+ H^s$ for $5/6 < s < 1$ , which do not have necessarily finite energy. The 
proof uses the I-method which is complicated by the fact that no 
$L^2$-conservation law holds. This shows that earlier results of Bethuel-Saut for data of the form $1 + H^1$ and 
G\'erard for finite energy data remain true for this class of rough data.
\end{abstract}
\maketitle
\renewcommand{\thefootnote}{\fnsymbol{footnote}}
\footnotetext{\hspace{-1.8em}{\it 2000 Mathematics Subject Classification:} 
35Q55, 35B60, 37L50}
\normalsize 
\setcounter{section}{0}
\section{Introduction and main results}
The Cauchy problem for the Gross-Pitaevskii 
equation in three space dimensions reads as follows 
\begin{eqnarray}
\label{0.1}
i\frac{\partial v}{\partial t} - \Delta v & = & v(1-|v|^2) 
\\
\label{0.2}
v(x,0) & = & v_0(x) \, ,
\end{eqnarray}
under the condition
\begin{equation}
\label{0.2a} 
 v \to 1  \quad {\mbox as} \quad |x| \to + \infty  \, , 
\end{equation}
where $v : {\mathbb R}^{1+3} \to {\mathbb C} $.

This problem occurs in theoretical physics, e.g.   
Bose-Einstein 
condensation and superfluidity, see \cite{Gr}, \cite{P},\cite{SS}.

Then one has the energy conservation law (see below)
\begin{equation}
\label{E0}
E(v(t)) = \int (|\nabla v(x,t)|^2 + \frac{1}{2} (|v(x,t)|^2 -1)^2) dx = E(v_0) 
\, .
\end{equation}
Because the solution does not vanish at infinity the standard theory for 
nonlinear Schr\"odinger equations is not directly applicable and it is natural 
to consider instead $u=v-1$ for a solution $v$ of (\ref{0.1}). Then $u$ 
satisfies the equivalent problem
\begin{eqnarray}
\label{0.3}
i\frac{\partial u}{\partial t} - \Delta u + (1+u)(|u|^2 + 2 \,Re\, u) & = & 0 
\\
\label{0.4}
u(x,0) & = & u_0(x) \, ,
\end{eqnarray}
under the condition
\begin{equation}
\label{0.4a} 
u \to 0  \quad {\mbox as} \quad |x| \to + \infty  \, . 
\end{equation}
The real part of the $L^2$-scalar product of equation (\ref{0.3}) with 
$\frac{\partial u}{\partial t}$ gives
$$ \frac{\partial}{\partial t} \int |\nabla u|^2 dx + 
\frac{1}{2}\frac{\partial}{\partial t} \int (|u(t)|^2 + 2 \, Re \, u(t))^2 dx 
= 
0 \, , $$
because 
$$ \frac{1}{2} \frac{\partial}{\partial t} ((|u|^2+2\, Re \, u)^2) = 2 \, Re 
\, 
((1+u)(|u|^2+2 \, Re \, u) \frac{\partial \bar{u}}{\partial t}) \, . $$
This gives the energy conservation law
$$ E(u(t)) = \int |\nabla u(t)|^2 dx + \frac{1}{2} \int (|u(t)|^2 + 2 \, Re \, 
u(t))^2 dx = E(u_0) \, . $$
In terms of $v$ one gets (\ref{E0}).
Remark that no conservation of $\|u(t)\|_{L^2}$ holds (in contrast to standard 
problems)!

We however get a bound for $\|u(t)\|_{L^2} = \|u(t)\|$ for finite energy data, which also belong to $L^2$, 
in the 
following way. The imaginary part of the scalar product of equation 
({\ref{0.3}}) with $u$ gives
$$ \frac{1}{2} \frac{\partial}{\partial t} \|u(t)\|^2 - \int (|u|^2 + 2 \, Re 
\, 
u) \, Im \, u \,  dx = 0 \, , $$
because
$$ Im \, (1+u)(|u|^2 + 2 \, Re \, u)\bar{u} = -(|u|^2+2 \, Re \, u) \, Im \, u 
\, . $$
This immediately implies
\begin{equation}
\label{**}
\frac{\partial}{\partial t} \|u(t)\|^2 \le 2 \int (|u(t)|^3 + 2|u(t)|^2) dx \, 
. 
\end{equation}
We also get
$$ \frac{\partial}{\partial t} \|u(t)\|^2 \le 2(\int(|u(t)|^2 + 2 \, Re \, 
u(t))^2 
dx)^{\frac{1}{2}} \|u(t)\| \le 2 \sqrt{2 E(u_0)} \|u(t)\| \, , $$
which implies
$$ \|u(t)\|^2 \le \|u_0\|^2 + \int_0^t 2 \sqrt{2E(u_0)} \|u(s)\| \, ds \, , $$
thus by a Gronwall type lemma
$$ \|u(t)\| \le \|u_0\| + \sqrt{2E(u_0)} t \, . $$
For data $u_0 \in H^1({\mathbb R}^3)$ these considerations lead directly to an 
a-priori-bound of $\|\nabla u(t)\|_{L^2}^2 \le E(u(t)) = E(u_0) $ , which is 
finite, because $H^1 \subset L^4$ by Sobolev's embedding theorem, and also to 
an 
a-priori bound of $\|u(t)\|_{L^2}$. Together with local well-posedness (cf. Theorem \ref{LWP} below) this shows that our problem 
(\ref{0.3}),(\ref{0.4}),(\ref{0.4a}) (and equivalently 
(\ref{0.1}),(\ref{0.2}),(\ref{0.2a})) has a unique global solution $u \in 
C^0({\mathbb R},H^1({\mathbb R}^3))$. 

The original proof was given by Bethuel and Saut \cite{BS}, Appendix A. Later 
G\'erard 
\cite{Ge} proved global 
well-posedness in the larger energy space using Strichartz estimates in two 
and three space dimensions. Gallo \cite{Ga} proved global well-posed- \\ness 
for more general nonlinearities for data with finite energy and space 
dimension $n \le 4$.

In the work at hand we are now interested in global well-posedness for data 
without finite energy, more precisely we consider solutions $v=1+u$, where $u 
\in H^s({\mathbb R}^3)$ for $s < 1$. We apply the so called I-method 
introduced 
by Colliander, Keel, Staffilani, Takaoka and Tao \cite{CKSTT} and successfully 
applied to various problems. There are two facts which 
complicate the problem: on one hand there is no scaling invariance and on the 
other hand no conservation law for the $L^2$-norm of $u$. As usual the energy 
conservation law is not directly applicable for $H^s$-data with $s<1$. However 
there is an "almost conservation law" for the modified energy $E(Iu)$, which 
is 
well defined for $u \in H^s$ (see the definition of $I$ below). This leads to 
an 
a-priori bound of $ \|\nabla Iu(t)\|_{L^2}$, if $s$ is close enough to 1, 
namely 
$s > 5/6$. This can be shown to be enough for an a-priori bound also for 
$\|u(t)\|_{L^2}$, which together gives a bound for $\|u(t)\|_{H^s}$. A local 
well-posedness result in Bourgain type spaces $X^{s,\frac{1}{2}+}[0,T] \subset 
C^0([0,T],H^s)$ with existence time dependent only on $\|u_0\|_{H^s}$ 
completes 
the global well-posedness result in this space. We even get unconditional 
global well-posedness in 
the space $C^0([0,T],H^s)$ using a result of Kato \cite{K}. This leads to the 
following 
main results (cf. the definition of the $X^{s,b}$-spaces below):
\begin{theorem}
\label{Theorem 1}
Let $ T > 0 $ , $ s > 5/6 $ and $u_0 \in H^s({\mathbb R}^3)$. The Cauchy problem 
(\ref{0.3}),(\ref{0.4}) has a unique global solution in 
$X^{s,\frac{1}{2}+}[0,T]$. This solution belongs to $C^0([0,T],H^s({\mathbb 
R}^3))$.
\end{theorem}
Combining this with the unconditional uniqueness result of T. Kato which we 
prove in Proposition \ref{Prop. 1} below we even get
\begin{theorem}
\label{Theorem 2}
Let $ T > 0$ , $ s > 5/6 $ and $u_0 \in H^s({\mathbb R}^3)$. The Cauchy problem 
(\ref{0.3}),\ref{0.4}) has a unique global solution in $C^0([0,T],H^s({\mathbb 
R}^3))$. Equivalently the Cauchy problem (\ref{0.1}),(\ref{0.2}) has a unique 
global solution in $C^0([0,T],1+H^s({\mathbb R}^3))$ for data $v_0 \in 
1+H^s({\mathbb R}^3)$.
\end{theorem}
The following proposition for more general nonlinearities and arbitrary 
dimensions goes back to Kato \cite{K}. We give the (short) proof in the 
special 
case of 
cubic polynomials as nonlinearity in three space dimensions.

\begin{prop}
\label{Prop. 1}
Assume $u_0 \in H^s({\mathbb R}^3)$. The Cauchy problem
$$ i \frac{\partial u}{\partial t} - \Delta u = F(u,\bar{u}) \quad , \quad 
u(0) 
= u_0 \, , $$
where $F(u,\bar{u})$ is a polynomial of degree three, has at most one solution 
$u \in C^0([0,T],H^s({\mathbb R}^3))$ for any $T>0$ , provided $ s \ge 2/3$ .
\end{prop}
\begin{proof}
Let $u,v \in C^0([0,T],H^s({\mathbb R}^3))$ be two solutions. By Sobolev's 
embedding 
$u,v \in C^0([0,T],L^{\frac{18}{5}})$ using $s \ge 2/3$. By the Strichartz 
estimates (see below) for the inhomogeneous Schr\"odinger equation we get 
(ignoring complex conjugates, which play no role here)
\begin{eqnarray*}
\lefteqn{ \|u-v\|_{L_t^3 L_x^{\frac{18}{5}}} + \|u-v\|_{L_t^{\infty} L_x^2} } 
\\
& \lesssim & \|u^3-v^3\|_{L_t^2 L_x^{\frac{6}{5}}} + \|u^2-v^2\|_{L_t^2 
L_x^{\frac{6}{5}}} + \|u-v\|_{L_t^1 L_x^2} \\
& \lesssim & \|u-v\|_{L_t^3 L_x^{\frac{18}{5}}} (\|u\|^2_{L^{12}_t 
L^{\frac{18}{5}}_x} + \|v\|^2_{L^{12}_t L^{\frac{18}{5}}_x}) \\
& & + \|u-v\|_{L_t^{\infty} L_x^2} (\|u\|_{L^2_t L^3_x} + \|v\|_{L^2_t L^3_x}) 
+ 
\|u-v\|_{L_t^1 L_x^2} \\
& \lesssim & \|u-v\|_{L_t^3 L_x^{\frac{18}{5}}} T^{\frac{1}{6}} 
(\|u\|^2_{L^{\infty}_t H^s_x} + \|v\|^2_{L^{\infty}_t H^s_x}) \\
& & + \|u-v\|_{L_t^{\infty} L_x^2} T^{\frac{1}{2}} (\|u\|_{L^{\infty}_t H^s_x} 
+ 
\|v\|_{L^{\infty}_t H^s_x} +1) \\
& \lesssim & \frac{1}{2} (\|u-v\|_{L_t^3 L_x^{\frac{18}{5}}} +  
\|u-v\|_{L_t^{\infty} L_x^2}) \, ,
\end{eqnarray*}
choosing $T$ small enough, which shows $u=v$.
\end{proof}

We use the following notation and well-known facts: the multiplier $I=I_N$ is 
for given $s<1$ and $N \ge 1$ defined by
$$ \widehat{I_N f}(\xi) := m_N(\xi) \widehat{f}(\xi) \, , $$
where $\, \widehat{}  \,$ denotes the Fourier transform with respect to the 
space 
variables. Here $m_N(\xi)$ is a smooth, radially symmetric, nonincreasing 
function of $|\xi|$ with 
$$ m_N(\xi) = \left\{              \begin{array}{ll}                   1 & 
|\xi| 
\le N\\   (\frac{N}{|\xi|})^{1-s}                 &    |\xi| \ge 2N           
\end{array}       \right. $$
We remark that $I: H^s \to H^1$ is a smoothing operator, so that especially 
$E(Iu)$ is well-defined for $u \in H^s({\mathbb R}^3)$ (remark that 
$H^1({\mathbb R}^3) \subset L^4({\mathbb R}^3)$).

We use the Bourgain type function space $X^{m,b}$ belonging to the 
Schr\"odinger 
equation $iu_t -\Delta u = 0$, which is defined as follows: let $\, \widehat{} 
\,$ or ${\mathcal F}$ denote 
the Fourier transform with respect to space and time and $\mathcal {F}^{-1}$ 
its 
inverse. $X^{m,b}$ is the 
completion of ${\mathcal S}({\mathbb R} \times {\mathbb R}^3)$ with respect to
$$ \|f\|_{X^{m,b}}  =  \| \langle \xi \rangle^m \langle \tau \rangle^b 
{\mathcal F}(e^{-it\Delta} f(x,t))\|_{L^2_{\xi, \tau}} = \| \langle \xi 
\rangle^m \langle \tau + |\xi|^2 \rangle^b 
\widehat{f}(\xi,\tau)\|_{L^2_{\xi,\tau}} \, , $$
For a given time interval $I$ we define
$$ \|f\|_{X^{m,b}(I)} := \inf_{g_{|I}=f}  \|g\|_{X^{m,b}} \, 
. $$
For $ s \ge 0$ and $1 \le r < \infty$ we denote by $H^{s,r}$ the standard Sobolev space, i.e. the completion of $C^{\infty}_0({\mathbb R}^3)$ with respect to
$$ \|f\|_{H^{s,r}} = \| {\mathcal F}^{-1}(\langle \xi \rangle^s \widehat{f}(\xi))\|_{L^r} \, .$$ 

We recall the following facts about the solutions u of the inhomogeneous 
linear 
Schr\"odinger equation (see e.g. \cite{GTV})
\begin{equation}
\label {4}
i u_t - \Delta u = F \quad , \quad u(0) = f \, . 
\end{equation}
For $ b'+1 \ge b \ge 0 \ge b' > -1/2 $ and 
$T \le 1$ we have
$$ \|u\|_{X^{m,b}[0,T]} \lesssim \|f\|_{H^m}  + 
T^{1+b'-b}\|F\|_{X^{m,b'}[0,T]} \, . $$
For $ 1/2 > b > b' \ge 0$ or $0 \ge b > b' > -1/2$: 
$$ \|f\|_{X^{m,b'}[0,T]} \lesssim T^{b-b'} \|f\|_{X^{m,b}[0,T]} $$
(see e.g. \cite{G}, Lemma 1.10).

Fundamental are the following Strichartz type estimates for the solution $u$ 
of 
(\ref{4}) in three space dimensions (see \cite{CH},\cite{KT}):
$$ \|u\|_{L^q(I,L^r({\mathbb R}^3))} \lesssim \|f\|_{L^2({\mathbb R}^3))} + 
\|F\|_{L^{\tilde{q}'}(I,L^{\tilde{r}'}({\mathbb R}^3))}$$
 with implicit constant independent of the interval $I \subset {\mathbb R}$
for all pairs $(q,r),(\tilde{q},\tilde{r})$ with $q,r,\tilde{q},\tilde{r} \ge 
2$ 
and $\frac{1}{q} + \frac{3}{2r} = \frac{3}{4}$ , $\frac{1}{\tilde{q}} + 
\frac{3}{2\tilde{r}} = \frac{3}{4}$, where 
$\frac{1}{\tilde{q}}+\frac{1}{\tilde{q}'} = 1$ and
$\frac{1}{\tilde{r}}+\frac{1}{\tilde{r}'} = 1$.
This implies
$$ \|\psi\|_{L^q(I,L^r({\mathbb R}^3))} \lesssim          
\|\psi\|_{X^{0,\frac{1}{2}+}(I)}  \, . $$
 
For real numbers $a$ we denote by $a+$, $a++$, $a-$ and $a--$ the numbers 
$a+\epsilon$, $a+2\epsilon$, $a-\epsilon$ and $a-2\epsilon$, respectively, 
where $\epsilon > 0$ is sufficiently small.

Of special interest is also a bilinear refinement, which goes back to 
Bourgain \cite{B}, namely the following frequency localized version in three 
dimensions:
\begin{lemma}
\label{Lemma}
Let $u_j$ be given with $supp \, \widehat{u}_j \subset \{|\xi| \sim N_j\}$ 
$(j=1,2)$ , $N_1 \le N_2$. Then the following estimates hold
\begin{eqnarray}
\label{B1}
\|u_1 u_2\|_{L^2_{x,t}}  \lesssim \frac{N_1}{N_2^{\frac{1}{2}}} 
\|u_1\|_{X^{0,\frac{1}{2}+}}\|u_2\|_{X^{0,\frac{1}{2}+}}  \, , \\
\label{B2}
\|u_1 u_2\|_{L^2_{x,t}}  \lesssim \frac{N_1^{1+}}{N_2^{\frac{1}{2}-}} 
\|u_1\|_{X^{0,\frac{1}{2}+}}\|u_2\|_{X^{0,\frac{1}{2}+}} \, .
\end{eqnarray} 
\end{lemma} 
\begin{proof}
For a proof of (\ref{B1}) we refer to Bourgain \cite{B}, Lemma 5 or Gr\"unrock 
\cite{G}. (\ref{B2}) follows 
by interpolation of (\ref{B1}) with the crude estimate
$$ \|u_1 u_2\|_{L^2_{xt}} \lesssim \|u_1\|_{L_t^{\infty} L_x^{6+}} 
\|u_2\|_{L^2_t L^{3-}_x} \lesssim N_1^{1+} 
\|u_1\|_{X^{0,\frac{1}{2}+}}\|u_2\|_{X^{0,\frac{1}{2}+}} \, , $$
using $X^{0,\frac{1}{2}+} \subset L_t^2 L_x^{3-}$ and Sobolev's embedding 
$\dot{H}^{1+} \subset L_x^{6+}$.
\end{proof}

The paper is organized as follows: in chapter 1 we prove two versions of a 
local 
well-posedness result for (\ref{0.3}),(\ref{0.4}), namely $u \in 
X^{s,\frac{1}{2}+}[0,\delta] $ for data $u_0 \in H^s$ with $s > 1/2$, and a 
modification 
where $\nabla Iu \in X^{0, \frac{1}{2}+}[0,\delta]$ for data $\nabla Iu_0 \in 
L^2$, which 
is necessary in order to combine it with an almost conservation law for the 
modified energy $E(Iu)$. In chapter 2 we use these local results and bounds 
for the modified energy given in chapter 3 in order to get the main theorem. 
It 
is namely shown that the bounds for the modified energy are enough to give 
also 
a uniform exponential bound for the $L^2$-norm of $u(t)$ and as a consequence 
for the $H^s$-norm for $u(t)$, which in view of the local well-posedness 
results 
suffices to get a global solution. In chapter 3 we calculate $\frac{d}{dt} 
E(Iu)$ for any solution of the equation $ i\frac{\partial Iu}{\partial t} - 
\Delta Iu + I((1+u)(|u|^2 + 2 \,Re\, u)) = 0 $. The most complicated part is 
to 
estimate the time integrated terms which appear in  $\frac{d}{dt} E(Iu)$. 
Finally we show that these estimates control the modified energy $E(Iu)$ 
uniformly on arbitrary time intervals $[0,T]$, provided $s > 5/6$.
\section{Local well-posedness}
The following local well-posedness theorem is more or less standard.
\begin{theorem}
\label{LWP}
Assume $s> 1/2$ and $u_0 \in H^s({\mathbb R}^3)$. Then the Cauchy problem 
(\ref{0.3}),(\ref{0.4}) is locally well-posed, i.e. there exists $T_0 = 
T_0(\|u_0\|_{H^s})$ such that there exists a unique solution $u \in 
X^{s,\frac{1}{2}+}(0,T_0)$. This solution belongs to $C^0([0,T_0],H^s({\mathbb 
R}^3))$. $T_0$ can be chosen such that $T_0 \sim 
\min(\|u_0\|_{H^s}^{-\frac{4}{2s-1}-},1)$.
\end{theorem}
\begin{proof}
We have to estimate $\|F(u)\|_{X^{s,-\frac{1}{2}+}}$ , where we define
\begin{equation}
\label{F}
F(u) = (1+u)(|u|^2 + 2 \,Re\, u)  \, .
\end{equation}
We want to show
$$ \||u|^2 u\|_{X^{s,-\frac{1}{2}++}} \lesssim T^{s-\frac{1}{2}-} 
\|u\|^3_{X^{s,\frac{1}{2}+}} \, , $$
where here and in the sequel we skip the interval $[0,T]$ in the 
$X^{s,b}[0,T]$-spaces. We ignore complex conjugates, because they play no role 
here, and use a fractional Leibniz rule and duality to reduce to the estimate
$$ \|u^2 \langle D \rangle^s u \psi \|_{L^1_{xt}} \lesssim T^{s-\frac{1}{2}-} 
\|u\|^3_{X^{s,\frac{1}{2}+}} \|\psi\|_{X^{0,\frac{1}{2}--}} \, . $$ 
We have
\begin{eqnarray*}
\|u^2  \langle D \rangle^s u \psi \|_{L^1_{xt}} 
& \lesssim & \|u^2 \langle D \rangle^s u \|_{L^{1+}_t L^2_x} 
\|\psi\|_{L^{\infty 
-}_t L^2_x} \\
& \lesssim & \|u\|^2_{L^{4+}_t L^6_x} \| \langle D \rangle^s u \|_{L^2_t 
L^6_x} 
\|\psi\|_{X^{0,\frac{1}{2}--}} \\
& \lesssim & \|u\|^2_{L^{4+}_t H^{s,r}_x} \|u\|_{X^{s,\frac{1}{2}+}} 
\|\psi\|_{X^{0,\frac{1}{2}--}} \\
& \lesssim & T^{s-\frac{1}{2}-} \|u\|^2_{L^{q}_t H^{s,r}_x} 
\|u\|_{X^{s,\frac{1}{2}+}} \|\psi\|_{X^{0,\frac{1}{2}--}} \\                     
  & \lesssim & T^{s-\frac{1}{2}-}  \|u\|_{X^{s,\frac{1}{2}+}}^3 
\|\psi\|_{X^{0,\frac{1}{2}--}} \, ,
\end{eqnarray*}                
where $\frac{1}{r} = \frac{1}{6} + \frac{s}{3}$ , so that $H^{s,r} \subset 
L^6$ 
, and $X^{0,\frac{1}{2}+} \subset L^2_t L^6_x$ by Strichartz, and $\frac{1}{q} 
= 
\frac{1-s}{2}$ , so that $\frac{2}{q} = 3(\frac{1}{2} - \frac{1}{r})$ , thus 
$X^{s,\frac{1}{2}+} \subset L^q_t H^{s,r}_x$ by Strichartz' estimate.

Similarly we get by Strichartz' estimate and Sobolev's embedding
$$
\|u \langle D \rangle^s u \|_{L^{1+}_t L^2_x}  \lesssim  T^{\frac{1}{2}-} \|u 
\langle D \rangle^s u \|_{L^2_t L^2_x} 
 \lesssim T^{\frac{1}{2}-} \|u\|_{L^{\infty}_t L^3_x} \| u\|_{L^2_t H^{s,6}_x} 
\lesssim T^{\frac{1}{2}-} \|u\|^2_{X^{s,\frac{1}{2}+}} \, , $$
thus
$$ \|u^2\|_{X^{s,-\frac{1}{2}++}} \lesssim T^{\frac{1}{2}-} 
\|u\|^2_{X^{s,\frac{1}{2}+}} \, . $$
Finally
$$ \|u\|_{X^{s,-\frac{1}{2}++}} \lesssim T^{\frac{1}{2}-} \|u\|_{L^2_t 
H^{s}_x} 
\lesssim T^{1-} \|u\|_{X^{s,\frac{1}{2}+}} \, . $$
Similar estimates hold for the difference 
$\|F(u)-F(v)|_{X^{0,-\frac{1}{2}+}}$. 
The standard Picard iteration shows the claimed result, where $ T \le 1 $ 
has to be chosen that $T^{s-\frac{1}{2}-} \|u_0\|_{H^s}^2 \lesssim 1$ and 
$T^{\frac{1}{2}-} \|u_0\|_{H^s} \lesssim 1$. Thus the choice as claimed in the 
theorem is possible.
\end{proof}
{\bf Remarks:}
1. A similar proof in spaces of the type $L^p_t L^r_x$ could also be given. 
This 
goes back to \cite{CW}, where $s=1/2$ is included, but in this limiting case 
the 
existence time depends not only on $\|u_0\|_{H^s}$.\\
2. Theorem \ref{LWP} shows that in order to get a global solution it is 
sufficient to have an a-priori bound of $\|u(t)\|_{H^s}$, if $s > 1/2$.

We next prove a similar local well-posedness result involving the operator 
$I$.
\begin{prop}
\label{LWP1}
Assume $ s > 1/2 $ and $\nabla Iu_0 \in L^2({\mathbb R}^3)$. Then (after 
application of $I$) the problem (\ref{0.3}),(\ref{0.4}) has a unique local 
solution $u$ with $\nabla Iu \in X^{0,\frac{1}{2}+}(0,\delta)$ and
$$\| \nabla Iu\|_{X^{0,\frac{1}{2}+}(0,\delta)} \le \sqrt{2} \|\nabla 
Iu_0\|_{L^2} \, , $$
where $ \delta \le 1 $ can be chosen such that
\begin{equation}
\label{Delta}
\left( \frac{\delta^{s-\frac{1}{2}-}}{N^{2(1-s)}} +  
\frac{\delta^{\frac{s}{2}-}}{N^{1-s}} + \delta^{\frac{1}{2}-} \right) \|\nabla 
Iu_0\|^2_{L^2} \sim 1 \, .
\end{equation}
\end{prop}
\begin{proof}
The cubic term in the nonlinearity will be estimated as follows:
\begin{equation}
\label{10}
\|\nabla I(u_1,u_2,u_3)\|_{X^{0,-\frac{1}{2}++}} \lesssim \left( 
\frac{\delta^{s-\frac{1}{2}-}}{N^{2(1-s)}} +  
\frac{\delta^{\frac{s}{2}-}}{N^{1-s}} + \delta^{\frac{1}{2}-} \right) \prod_{i=1}^3   \|\nabla Iu_i\|_{X^{0,\frac{1}{2}+}}
\end{equation} 
This follows from
\begin{eqnarray*}
A & := & \int_0^{\delta}\int_* M(\xi_1,\xi_2,\xi_3) \prod_{i=1}^3 
\widehat{u_i}(\xi_i,t) \widehat{\psi}(\xi_4,t) d\xi_1 d\xi_2 d\xi_3 d\xi_4 dt 
\\
& & \lesssim \left( \frac{\delta^{s-\frac{1}{2}-}}{N^{2(1-s)}} +  
\frac{\delta^{\frac{s}{2}-}}{N^{1-s}} + \delta^{\frac{1}{2}-} \right) 
\prod_{i=1}^3 \|u_i\|_{X^{0,\frac{1}{2}+}} \|\psi\|_{X^{0,\frac{1}{2}--}} \,
\end{eqnarray*}
where
$$
M(\xi_1,\xi_2,\xi_3) := \frac{m(\xi_1+\xi_2+\xi_3)}{m(\xi_1)m(\xi_2)m(\xi_3)} 
\frac{|\xi_1+\xi_2+\xi_3|}{|\xi_1||\xi_2||\xi_3|}
$$
and * denotes integration over the region $\{\sum_{i=1}^4 \xi_i = 0\}$. We 
assume here and in the following w.l.o.g. that the Fourier transforms are 
nonnnegative. 
We also assume w.l.o.g. $|\xi_1| \ge |\xi_2| \ge |\xi_3|$. \\
{\bf Case 1:} $|\xi_1| \ge |\xi_2| \ge |\xi_3| \ge N$. \\
We first estimate the multiplier $M$. 
If $|\xi_1+\xi_2+\xi_3| \ge N$ we get
\begin{eqnarray*}
M(\xi_1,\xi_2,\xi_3) & \lesssim &\prod_{i=1}^3 (\frac{|\xi_i|}{N})^{1-s} 
\frac{N^{1-s}}{|\xi_1 + \xi_2 +\xi_3|^{1-s}} 
\frac{|\xi_1+\xi_2+\xi_3|}{|\xi_1||\xi_2||\xi_3|} \\
& \lesssim & \prod_{i=1}^3 (\frac{|\xi_i|}{N})^{1-s} \frac{N^{1-s}|\ 
\xi_1|^s}{|\xi_1||\xi_2||\xi_3|} \lesssim \frac{1}{|\xi_2|^s |\xi_3|^s 
N^{2(1-s)}} \, ,
\end{eqnarray*}
and if $|\xi_1+\xi_2+\xi_3| \le N$ we have
$$ M(\xi_1,\xi_2,\xi_3) \lesssim \prod_{i=1}^3 (\frac{|\xi_i|}{N})^{1-s} 
\frac{N}{|\xi_1||\xi_2||\xi_3|} \lesssim\prod_{i=1}^3 \frac{N^s}{|\xi_i|^s 
N^{2(1-s)}} \lesssim \frac{1}{|\xi_2|^s |\xi_3|^s N^{2(1-s)}}$$
as before. This implies by H\"older's and Strichartz' inequality and Sobolev's 
embedding
\begin{eqnarray*}
A & \lesssim & \frac{1}{N^{2(1-s)}} \|\psi\|_{L_t^{\tilde{q}} L_x^2}
\|u_1\|_{L_t^2 L_x^6} \|{\mathcal 
F}^{-1}(\frac{\widehat{u_2}}{|\xi_2|^s})\|_{L_t^q L_x^6} 
\|{\mathcal F}^{-1}(\frac{\widehat{u_3}}{|\xi_3|^s})\|_{L_t^q L_x^6} \\
& \lesssim & \frac{\delta^{s-\frac{1}{2}-}}{N^{2(1-s)}} \|\psi\|_{L_t^{\infty 
-} 
L_x^2} \|u_1\|_{X^{0,\frac{1}{2}+}} \|u_2\|_{L^q_t L_x^r}  \|u_3\|_{L^q_t 
L_x^r} 
\\
& \lesssim & \frac{\delta^{s-\frac{1}{2}-}}{N^{2(1-s)}} 
\|\psi\|_{X^{0,\frac{1}{2}--}} \prod_{i=1}^3 \|u_i\|_{X^{0,\frac{1}{2}+}} \, .
\end{eqnarray*}    
where $\frac{1}{q} = \frac{1}{2}-\frac{s}{2}$ , $\frac{1}{\tilde{q}} = 
s-\frac{1}{2}$, such that $\dot{H}^{s,r} \subset L^6$ with 
$\frac{1}{r}=\frac{1}{6}+\frac{s}{3}$ and $\frac{1}{q} = 
\frac{3}{4}-\frac{3}{2r}$.\\
{\bf Case 2:} $|\xi_1| \ge |\xi_2| \ge N \ge |\xi_3|$.\\
The multiplier $M$ is estimated as follows: if $|\xi_1+\xi_2+\xi_3| \ge N$ we 
get
$$M(\xi_1,\xi_2,\xi_3) \lesssim  (\frac{|\xi_1|}{N})^{1-s} 
(\frac{|\xi_2|}{N})^{1-s} \frac{N^{1-s}}{|\xi_1 + \xi_2 +\xi_3|^{1-s}} 
\frac{|\xi_1+\xi_2+\xi_3|}{|\xi_1||\xi_2||\xi_3|} 
 \lesssim \frac{1}{|\xi_2|^s |\xi_3| N^{1-s}} \, ,$$
and if $|\xi_1+\xi_2+\xi_3| \le N$ we also have
$$ M(\xi_1,\xi_2,\xi_3) \lesssim (\frac{|\xi_1|}{N})^{1-s} 
(\frac{|\xi_2|}{N})^{1-s} \frac{N}{|\xi_1||\xi_2||\xi_3|} \lesssim 
\frac{1}{|\xi_2|^s |\xi_3| N^{1-s}} \, .$$
This implies by H\"older, Strichartz and Sobolev
\begin{eqnarray*}
A & \lesssim & \frac{1}{N^{1-s}} \|\psi\|_{L_t^{\frac{2}{s}} L_x^2} 
\|u_1\|_{L_t^2 L_x^6} \|{\mathcal F}^{-1}(|\xi_2|^{-s} \widehat{u_2})\|_{L_t^q 
L_x^6} \|{\mathcal F}^{-1}(|\xi_3|^{-1} \widehat{u_3})\|_{L_t^{\infty} L_x^6} 
\\
& \lesssim & \frac{1}{N^{1-s}} \delta^{\frac{s}{2}-} \|\psi\|_{L_t^{\infty-} 
L_x^2} \|u_1\|_{L_t^2 L_x^6} \|u_2\|_{L_t^q L_x^r} \|u_3\|_{L_t^{\infty} 
L_x^2} 
\\
& \lesssim & \frac{1}{N^{1-s}} \delta^{\frac{s}{2}-} 
\|\psi\|_{X^{0,\frac{1}{2}--}} \prod_{i=1}^3 \|u_i\|_{X^{0,\frac{1}{2}+}} \, ,
\end{eqnarray*}
where $\frac{1}{q} = \frac{1}{2}-\frac{s}{2}$, $\dot{H}^{s,r}_x \subset L^6_x$ 
for $\frac{1}{r} = \frac{1}{6}+\frac{s}{3}$, thus $\frac{1}{q} = 
\frac{3}{4}-\frac{3}{2r}$ and $X^{0,\frac{1}{2}+} \subset L_t^q L_x^r$.\\
{\bf Case 3:} $|\xi_1| \ge N \ge |\xi_2| \ge |\xi_3|$ and $|\xi_1| \gg 
|\xi_2|$ 
, or $N \ge |\xi_1| \ge |\xi_2| \ge |\xi_3|$. \\
In these cases we have $M(\xi_1,\xi_2,\xi_3) \lesssim \frac{1}{|\xi_2| 
|\xi_3|}$, thus
\begin{eqnarray*}
A & \lesssim &  \|\psi\|_{L^2_{xt}} \|u_1\|_{L_t^2 L_x^6} \|{\mathcal 
F}^{-1}(|\xi_2|^{-1} \widehat{u_2})\|_{L_t^{\infty} L_x^6} \|{\mathcal 
F}^{-1}(|\xi_3|^{-1} \widehat{u_3})\|_{L_t^{\infty} L_x^6} \\
& \lesssim & \delta^{\frac{1}{2}-} \|\psi\|_{X^{0,\frac{1}{2}--}} 
\prod_{i=1}^3 
\|u_i\|_{X^{0,\frac{1}{2}+}} 
\end{eqnarray*}
similarly as in case 2. This implies (\ref{10}).

Next we have to estimate the quadratic terms in the nonlinearity. We want to 
show
\begin{equation}
\label{11}
\|\nabla I(u_1 u_2)\|_{X^{0,-\frac{1}{2}++}} \lesssim \delta^{\frac{1}{2}-} 
\|\nabla Iu_1\|_{X^{0,\frac{1}{2}+}} \|\nabla Iu_2\|_{X^{0,\frac{1}{2}+}} \, , 
\end{equation}
which follows from
\begin{eqnarray*}
B & := & \int_0^{\delta} \int_* M(\xi_1,\xi_2) \widehat{u_1}(\xi_1,t) 
\widehat{u_2}(\xi_1,t) \widehat{\psi}(\xi_3,t) d\xi_1 d\xi_2 d\xi_3 dt \\
& \lesssim & \delta^{\frac{1}{2}-} \|u_1\|_{X^{0,\frac{1}{2}+}}
\|u_2\|_{X^{0,\frac{1}{2}+}} \|\psi\|_{X^{0,\frac{1}{2}--}} \, , 
\end{eqnarray*}
where * denotes integration over the region $\{\xi_1 + \xi_2 + \xi_3 = 0 \}$
and
$$ M(\xi_1,\xi_2) := \frac{m(\xi_1 + \xi_2)}{m(\xi_1) m(\xi_2)} \frac{|\xi_1 + 
\xi_2|}{|\xi_1||\xi_2|} \, . $$
Assuming w.l.o.g. $|\xi_1| \ge |\xi_2|$ we first consider\\
{\bf Case 1:} $|\xi_2| \ge N$.\\
 If $|\xi_1 + \xi_2| \ge N$ we have
$$ M(\xi_1,\xi_2) \lesssim (\frac{|\xi_1|}{N})^{1-s} (\frac{|\xi_2|}{N})^{1-s}
\frac{N^{1-s}}{|\xi_1 + \xi_2|^{1-s}} \frac{|\xi_1 + \xi_2|}{|\xi_1||\xi_2|} 
\lesssim 
\frac{1}{|\xi_2|^s N^{1-s}} \, ,$$
and in the case $|\xi_1 + \xi_2| \le  N$ we also get
$$  M(\xi_1,\xi_2) \lesssim (\frac{|\xi_1|}{N})^{1-s} 
(\frac{|\xi_2|}{N})^{1-s}
\frac{N}{|\xi_1||\xi_2|} \lesssim \frac{1}{N^{2(1-s)}} \frac{N^s 
N^{1-s}}{|\xi_2|^s |\xi_1|^s} 
\lesssim 
\frac{1}{|\xi_2|^s N^{1-s}} \, ,$$
so that by Strichartz' estimate using $s \ge\frac{1}{2}$:
\begin{eqnarray*}
B & \lesssim & \frac{1}{N^{1-s}} \|\psi\|_{L^2_{xt}} \|u_1\|_{L_t^2 L_x^6} 
\|{\mathcal F}^{-1}(\frac{\widehat{u_2}}{|\xi_2|^s})\|_{L^{\infty}_t L^3_x} \\
& \lesssim & \frac{\delta^{\frac{1}{2}-}}{N^{1-s}} 
\|\psi\|_{X^{0,\frac{1}{2}--}}
 \|u_1\|_{X^{0,\frac{1}{2}+}} \|{\mathcal 
F}^{-1}(\frac{\widehat{u_2}}{|\xi_2|^s})\|_{L^{\infty}_t
 \dot{H}^{\frac{1}{2}}_x} \\ 
& \lesssim &\frac{\delta^{\frac{1}{2}-}}{N^{1-s}} 
\|\psi\|_{X^{0,\frac{1}{2}--}}
 \|u_1\|_{X^{0,\frac{1}{2}+}} \|u_2\|_{X^{0,\frac{1}{2}+}} \, .
\end{eqnarray*}
{\bf Case 2:} $N \ge |\xi_2|$. \\
If $|\xi_1| \gg |\xi_2|$ we have $ M(\xi_1+\xi_2) \sim 
\frac{m(\xi_1 + \xi_2)}{m(\xi_1)m(\xi_2)|\xi_2|} \sim \frac{1}{|\xi_2|} \, , $
whereas, if $|\xi_1| \sim |\xi_2|$, we have $|\xi_1 +\xi_2| \lesssim N$ and 
$m(\xi_1) \sim m(\xi_2) \sim m(\xi_1 + \xi_2) \sim 1$, which leads to the same 
bound for $M(\xi_1,\xi_2)$. Thus
$$ B \lesssim \|\psi\|_{L^2_{xt}} \|u_1\|_{L^2_t L_x^3}\|{\mathcal 
F}^{-1}(\frac{\widehat{u_2}(\xi_2)}{|\xi_2|})\|_{L^{\infty}_t L_x^6} \lesssim 
\delta^{\frac{1}{2}-} \|\psi\|_{X^{0,\frac{1}{2}--}} 
\|u_1\|_{X^{0,\frac{1}{2}+}} \|u_2\|_{X^{0,\frac{1}{2}+}} \, . $$

Finally 
$$ \|\nabla Iu\|_{X^{0,-\frac{1}{2}++}} \lesssim \delta^{\frac{1}{2}-} 
\|\nabla 
Iu\|_{X^{0,0}} \lesssim \delta^{1-} \|\nabla Iu\|_{X^{0,\frac{1}{2}+}} \, . $$

Similar estimates hold for the difference $ \| \nabla 
I(F(u)-F(v))\|_{X^{0,-\frac{1}{2}+}} $.

A Picard iteration leads to the desired solution in $[0,\delta]$, where 
$\delta 
\le 1$ has to be chosen such that
$$ \frac{\delta^{s-\frac{1}{2}-}}{N^{2(1-s)}} \|\nabla Iu_0\|_{L^2}^2 \lesssim 
1 
\quad , \quad \delta^{\frac{1}{2}-} \|\nabla Iu_0\|_{L^2}^2 \lesssim 1 \, ,$$ 
$$ 
 \frac{\delta^{\frac{s}{2}-}}{N^{1-s}} \|\nabla Iu_0\|_{L^2}^2 \lesssim 1 
\quad 
, \quad \delta^{\frac{1}{2}-} \|\nabla Iu_0\|_{L^2} \lesssim 1 \, . $$
\end{proof}
{\bf Remark:} We want to iterate this local existence theorem with time steps 
of 
equal length until we reach a given (large) time $T$. For this we need to 
control
\begin{equation}
\label{6}
\|\nabla Iu(t)\|_{L^2} \le c(T)   \quad \forall \, 0 \le t \le T \, . 
\end{equation}
This is achieved for $u_0 \in H^s$ and $s > 5/6$ by giving uniform bounds of 
the 
modified energy $E(Iu(t))$, which is done in chapter 3.

\section{ Proof of Theorem \ref{Theorem 1}}
\begin{proof}
Let us assume for the moment that (\ref{6}) holds and show that this leads to 
the claim of Theorem \ref{Theorem 1}. We thus have an a-priori bound for our 
local solution of Proposition \ref{LWP1} on any existence interval $[0,T]$, 
namely of
 \begin{equation}
\label{7} 
\|\nabla Iu(t)\|_{L^2} \sim \||\xi| \widehat{u}(\xi,t)\|_{L^2(\{|\xi| \le 
N\})} 
+ \||\xi|^s \widehat{u}(\xi,t)\|_{L^2(\{|\xi| \ge N\})} N^{1-s} \, . 
\end{equation}
What remains to be given is an a-priori bound for $\|u(t)\|_{L^2}$ as a 
consequence of (\ref{**}).
\begin{lemma}
\label{lemma}
On any existence interval $[0,T]$ of our solution $u \in 
X^{s,\frac{1}{2}+}[0,T]$ we have $ \|u(t)\|_{L^2({\mathbb R}^3)} \le c(T)$ .
\end{lemma}
\begin{proof}
We smoothly decompose $\widehat{u} = \widehat{u_1} + \widehat{u_2}$ with $ 
supp 
\, \widehat{u_1} \subset \{ |\xi| \le 2 \}$ and $ supp \, \widehat{u_2} 
\subset 
\{|\xi| \ge 1\}$. Then we have by Gagliardo-Nirenberg
\begin{eqnarray*}
\|u\|_{L^3} & \le &  \|u_1\|_{L^3} + \|u_2\|_{L^3} \lesssim \|\nabla 
u_1\|_{L^2}^{\frac{1}{2}} \|u_1\|_{L^2}^{\frac{1}{2}} + \| |D|^{\frac{1}{2}} 
u_2\|_{L^2} \\
& \lesssim &  \|\nabla u_1\|_{L^2}^{\frac{1}{2}} \|u_1\|_{L^2}^{\frac{1}{2}} + 
\| |D|^s u_2\|_{L^2}^{\frac{1}{2s}} \|u_2\|_{L^2}^{1-\frac{1}{2s}} \\
& \lesssim &  \|\nabla u_1\|_{L^2}^2 + \|u_1\|_{L^2}^{\frac{2}{3}} +  
\|u_2\|_{L^2}^{\frac{2}{3}}+ \| |D|^s u_2\|_{L^2}^{\frac{2}{3-2s}} \, ,
\end{eqnarray*}
so that by (\ref{6}) and (\ref{7}) we get on $[0,T]$:
\begin{eqnarray*}
\|u(t)\|_{L^3}^3 &  \lesssim & \||\xi| \widehat{u_1}(\xi,t)\|_{L^2}^6 +
\|u_1(t)\|_{L^2}^2 + \|u_2(t)\|_{L^2}^2 + \||\xi|^s 
\widehat{u_2}(\xi,t)\|_{L^2}^{\frac{6}{3-2s}}\\
& \le &  c'(T)( \|u(t)\|_{L^2}^2 +1) \, . 
\end{eqnarray*}
(\ref{**}) gives
$$ \frac{d}{dt} \|u(t)\|_{L^2}^2 \le c''(T)( \|u(t)\|_{L^2}^2 +1) \, , 
$$
so that Gronwall's lemma gives
$$ \|u(t)\|_{L^2}^2 + 1  \le (\|u_0\|_{L^2}^2 +1) e^{c''(T)T} $$
on $[0,T]$.
\end{proof}
Combining Lemma \ref{lemma} with  (\ref{6}) and (\ref{7})  we get an a-priori 
bound of 
$\|u(t)\|_{H^s}$. Together with Theorem \ref{LWP} we immediately get Theorem 
\ref{Theorem 1}.
\end{proof}

\section{Estimates for the modified energy}
Application of the operator $I$ to equation (\ref{0.3}) gives
\begin{equation}
\label{5}
i \frac{\partial}{\partial t} Iu - \Delta Iu + IF(u) = 0 \, ,
\end{equation}
with
$$F(u) := (1+u)(|u|^2 + 2 \, Re \, u) \, . $$
We define the modified energy 
$$ E(Iu) = \int |\nabla Iu|^2 dx + \frac{1}{2} \int (|Iu|^2 + 2 \, Re \, Iu)^2 
dx \, . $$
Of course one cannot expect that it is conserved, but we want to show an 
almost 
conservation law for it. We calculate its derivative as follows:
$$ \frac{d}{dt} E(Iu) = 2 \, Re \, \langle -\Delta Iu + (|Iu|^2+ 2\, Re\, 
Iu)(1+Iu),Iu_t \rangle = 2 \, Re \, \langle F(Iu)-IF(u),Iu_t \rangle $$
by replacing $\Delta Iu$ using (\ref{5}). Next we replace $Iu_t$ again by use 
of 
(\ref{5}) and get
\begin{eqnarray}
\nonumber
\frac{d}{dt} E(Iu)& =& Im \, (\langle \nabla(F(Iu)-IF(u)),\nabla Iu \rangle + 
\langle F(Iu)-IF(u),IF(u) \rangle ) \\
\label{E}
& \le & |\langle \nabla(F(Iu)-IF(u)),\nabla Iu \rangle| + |\langle 
F(Iu)-IF(u),IF(u) \rangle | \, .
\end{eqnarray}
In order to control the increment of $E(Iu)$ by (\ref{E}) on the local 
existence 
interval $[0,\delta]$ we have to estimate several terms. We assume from now on 
$ 
s \ge 3/4$.\\
{\bf 1.} Let us first consider the first term on the right hand side
of (\ref{E}). Here and in the following we ignore complex conjugates,
because they are of no interest. We want to show
$$ \int_0^{\delta} |\langle \nabla(I(u^3)-(Iu)^3),\nabla Iu \rangle| dt 
\lesssim 
N^{-1+} \|\nabla Iu\|^4_{X^{0,\frac{1}{2}+}} \, . $$
This follows from
\begin{equation}
\label{12}
A = | \int_0^{\delta} \int_* M(\xi_1,\xi_2,\xi_3) \prod_{i=1}^4 
\widehat{u_i}(\xi_i,t) d\xi_1 d\xi_2 d\xi_3 d\xi_4 dt | \lesssim N^{-1+} 
\prod_{i=1}^4 \|u_i\|_{X^{0,\frac{1}{2}+}} \, , 
\end{equation}
where * denotes integration over the region $\{\sum_{i=1}^4 \xi_i = 0 \}$ and
$$ M(\xi_1,\xi_2,\xi_3) := \frac{|m(\xi_1 + \xi_2 + \xi_3) - m(\xi_1) m(\xi_2) 
m(\xi_3)|}{ m(\xi_1) m(\xi_2) m(\xi_3)} \frac{|\xi_1 + \xi_2 + 
\xi_3|}{|\xi_1||\xi_2||\xi_3|} \, . $$
We assume here and in the following that the Fourier transforms are 
nonnegative w.l.o.g.
In most of the cases we perform dyadic decompositions with respect to 
$|\xi_i|$, 
where $|\xi_i|\sim N_i$ with $N_i = 2^{k_i}$ , $k_i \in {\mathbb Z}$. In order 
to sum the dyadic parts at the end we always need a convergence generating 
factor $\frac{1 \wedge N_{min}^{0+}}{N_{max}^{0+}}$, where $N_{min}$ and 
$N_{max}$ is the smallest and the largest of the numbers $N_i$, respectively. 
$N_{max} \ge N \ge 1$  can be assumed in all the cases, because otherwise our 
multiplier $M$ is identically zero.

In the term at hand we also assume w.l.o.g. $N_1 \ge N_2 \ge N_3$ and $N_1 \ge 
N$.\\
{\bf Case 1:} $N_1 \ge N_2 \ge N_3 \gtrsim N$ $\Rightarrow N_4 \lesssim N_1 
\sim 
N_{max}$. \\
We have for $s\ge 3/4$:
$$M(\xi_1,\xi_2,\xi_3) \lesssim (\frac{N_1}{N})^{\frac{1}{4}}  
(\frac{N_2}{N})^{\frac{1}{4}}  (\frac{N_3}{N})^{\frac{1}{4}} \frac{N_1}{N_1 
N_2 
N_3} \, . $$
Thus by the bilinear Strichartz estimate (\ref{B2}) we get
\begin{eqnarray*}
A & \lesssim &  (\frac{N_1}{N})^{\frac{1}{4}}  (\frac{N_2}{N})^{\frac{1}{4}}  
(\frac{N_3}{N})^{\frac{1}{4}} \frac{1}{ N_2 N_3} \|u_1 u_3\|_{L^2_{xt}} 
\|u_2\|_{L^{\infty}_t L^{3-}_{x}} \|u_4\|_{L^2_t L^{6+}_x} \\
 & \lesssim &  (\frac{N_1}{N})^{\frac{1}{4}}  (\frac{N_2}{N})^{\frac{1}{4}}  
(\frac{N_3}{N})^{\frac{1}{4}} \frac{1}{ N_2 N_3} 
\frac{N^{1+}_3}{N_1^{\frac{1}{2}-}}  N_2^{\frac{1}{2}-}  N_4^{0+}   
\prod_{i=1}^4 \|u_i\|_{X^{0,\frac{1}{2}+}} \\
& \lesssim & \frac{N_3^{\frac{1}{4}+} N_4^{0+}}{N_1^{\frac{1}{4}-} 
N_2^{\frac{1}{4}+} N^{\frac{3}{4}}} \prod_{i=1}^4 \|u_i\|_{X^{0,\frac{1}{2}+}} 
\\
& \lesssim & \frac{1 \wedge N_{min}^{0+}}{ N_{max}^{0+} N^{1-}} \prod_{i=1}^4 
\|u_i\|_{X^{0,\frac{1}{2}+}} \, .
\end{eqnarray*}
{\bf Case 2:} $N_1 \ge N_2 \gtrsim N \gtrsim N_3$. \\
This gives the same bound as in case 1 (without the factor 
$(\frac{N_3}{N})^{\frac{1}{4}}$). \\
{\bf Case 3:} $N_1 \gtrsim N \gtrsim N_2 \ge N_3$ and $ N_1 \gg N_2$. \\
By the mean value theorem we get
$$M(\xi_1,\xi_2,\xi_3) \lesssim \frac{|\nabla m(\xi_1) \xi_2|}{m(\xi_1)} \frac{N_1}{N_1 N_2 N_3}
\lesssim \frac{N_2}{N_1} \frac{N_1}{N_1 N_2 N_3}  $$
leading as in case 1 to the bound
\begin{eqnarray*}
A &  \lesssim & \frac{N_2}{N_1} \frac{N_1}{N_1 N_2 N_3} \frac{N_3^{1+} 
N_2^{\frac{1}{2}-} N_4^{0+}}{N_1^{\frac{1}{2}-}}  \prod_{i=1}^4 
\|u_i\|_{X^{0,\frac{1}{2}+}} \\
& \lesssim & \frac{N_2^{\frac{1}{2}-}N_3^{0+} N_4^{0+}}{N_1^{\frac{3}{2}-}}  
\prod_{i=1}^4 \|u_i\|_{X^{0,\frac{1}{2}+}} \lesssim \frac{1 \wedge 
N_{min}^{0+}}{ N_{max}^{0+} N^{1-}} \prod_{i=1}^4 \|u_i\|_{X^{0,\frac{1}{2}+}} 
\, .
\end{eqnarray*}
This proves (\ref{12}) after dyadic summation over $N_1,N_2,N_3,N_4$.\\
{\bf 2.} We next want to show
$$ \int_0^{\delta} |\langle \nabla(I(u^2)-(Iu)^2),\nabla Iu \rangle| dt 
\lesssim 
N^{-1+} \delta^{\frac{1}{4}-} \|\nabla Iu\|^3_{X^{0,\frac{1}{2}+}} \, $$
which follows from
\begin{equation}
\label{13}
B = | \int_0^{\delta} \int_* M(\xi_1,\xi_2,\xi_3) \prod_{i=1}^3 
\widehat{u_i}(\xi_i,t) d\xi_1 d\xi_2 d\xi_3 dt | \lesssim 
N^{-1+}\delta^{\frac{1}{4}-} \prod_{i=1}^3 \|u_i\|_{X^{0,\frac{1}{2}+}} \, , 
\end{equation}
where * denotes integration over the region $\{\sum_{i=1}^3 \xi_i = 0 \}$ and
$$ M(\xi_1,\xi_2,\xi_3) := \frac{|m( \xi_2 + \xi_3) -  m(\xi_2) m(\xi_3)|}{ 
m(\xi_2) m(\xi_3)} \frac{ |\xi_2 + \xi_3|}{|\xi_2||\xi_3|} \, . $$
Assume w.l.o.g. $|\xi_2| \ge |\xi_3|$ and $|\xi_2| \ge N$ $\Rightarrow |\xi_1| 
\lesssim |\xi_2| \sim N_{max}$.\\
{\bf Case 1:} $N_2 \sim N_3 \gtrsim N$.\\
We have
$$ M(\xi_1,\xi_2,\xi_3) \lesssim \frac{1}{|m(\xi_2)|^2 N_2} \lesssim 
(\frac{N_2}{N})^{\frac{1}{2}} \frac{1}{N_2} \, . $$
Using $\dot{H}^{0+,3}_x \subset L_x^{3+}$ and $X^{0,\frac{1}{2}+} \subset 
L_t^{4+} L_x^{3-}$ we get
\begin{eqnarray*}
B & \lesssim & (\frac{N_2}{N})^{\frac{1}{2}} \frac{1}{N_2} 
\|u_1\|_{L^{3+}_{xt}} 
\|u_2\|_{L^{3-}_{xt}} \|u_3\|_{L^{3+}_{xt}} 
\lesssim\frac{\delta^{\frac{1}{4}-} 
N_1^{0+} N_3^{0+}}{N_{max}^{0+} N^{1-}} \prod_{i=1}^3 \|u_i\|_{L^{4+}_t 
L^{3-}_x} \\
& \lesssim & \frac{\delta^{\frac{1}{4}-}(1 \wedge N_{min}^{0+}) }{N_{max}^{0+} 
N^{1-}} \prod_{i=1}^3 \|u_i\|_{X^{0,\frac{1}{2}+}} \, .
\end{eqnarray*}
{\bf Case 2:} $N_1 \sim N_3$ $\Rightarrow N_2 \lesssim N_1$.\\
Because (see above) $N_1 \lesssim N_2$ we have $N_1 \sim N_2 \sim N_3$, which 
was already considered in case 1.\\
{\bf Case 3:} $|\xi_1| \sim |\xi_2|$ $\Rightarrow |\xi_3| \lesssim |\xi_1|$ 
and 
$|\xi_3| \ll |\xi_1| \sim |\xi_2|$.\\
a. If $|\xi_3| \gtrsim N$, we get as in case 1
$$M(\xi_1,\xi_2,\xi_3) \lesssim \frac{1}{m(\xi_3) |\xi_3|} \lesssim 
(\frac{|\xi_3|}{N})^{\frac{1}{4}} \frac{1}{|\xi_3|} \lesssim \frac{1}{N} \, . 
$$
Thus by Strichartz' estimate
$$ B \lesssim \frac{1}{N} \prod_{i=1}^3 \|u_i\|_{L^3_{xt}} \lesssim 
\frac{\delta^{\frac{1}{4}}}{N} \prod_{i=1}^3 \|u_i\|_{L^4_t L^3_x} \lesssim 
\lesssim \frac{\delta^{\frac{1}{4}}}{N} \prod_{i=1}^3 
\|u_i\|_{X^{0,\frac{1}{2}+}} \, . $$
b. If $|\xi_3| \lesssim N$ we estimate
$$ M(\xi_1,\xi_2,\xi_3) \sim \frac{|(\nabla m)(\xi_2) \xi_3|}{m(\xi_2) 
|\xi_3|} 
\sim \frac{1}{|\xi_2|} \lesssim \frac{1}{N} $$
leading to the same bound as in a. \\
This proves (\ref{13}).\\
{\bf 3.} Next we consider the second term on the right hand side of (\ref{E}) 
and want to show
$$ \int_0^{\delta} |\langle I(u^3) - (Iu)^3,I(u^3)\rangle| dt \lesssim N^{-2+} 
\|\nabla Iu\|_{X^{0,\frac{1}{2}+}}^6 \, . $$
This means that we have to show
\begin{eqnarray}
\label{14}
C  =  |\int_0^{\delta} \int_* M(\xi_1,...,\xi_6) \prod_{i=1}^6 
\widehat{u_i}(\xi_i,t)d\xi_1 ...d\xi_6 dt| \lesssim N^{-2+} \prod_{i=1}^6 
\|u_i\|_{X^{0,\frac{1}{2}+}} \, , 
\end{eqnarray}
where
$$ M(\xi_1,...,\xi_6) := 
\frac{|m(\xi_1+\xi_2+\xi_3)-m(\xi_1)m(\xi_2)m(\xi_3)|}{m(\xi_1)m(\xi_2)m(\xi_3)} 
\frac{m(\xi_4+\xi_5+\xi_6)}{m(\xi_4)m(\xi_5)m(\xi_6)} \prod_{i=1}^6 
|\xi_i|^{-1} 
\, . $$
Assume w.l.o.g. $N_1\ge  N_2 \ge N_3$ , $N_1 \ge N$ and $N_4 \ge N_5 \ge 
N_6$.\\
{\bf Case 1:} $N_4 \ge N_5 \ge N_6 \ge N$.\\
{\bf a.} $N_1 \ge N_2 \ge N_3 \ge N$.\\
We have
$$ M(\xi_1,...,\xi_6) \lesssim \prod_{i=1}^6 (\frac{N_i}{N})^{\frac{1}{4}} 
\prod_{i=1}^6 N_i^{-1} \,. $$
Thus by Strichartz and Sobolev:
\begin{eqnarray*}
C & \lesssim &  \prod_{i=1}^6 (\frac{N_i}{N})^{\frac{1}{4}} \prod_{i=1}^6 
N_i^{-1} \|u_1\|_{L_t^2 L_x^6} \|u_2\|_{L_t^{\infty} 
L_x^3}\|u_4\|_{L_t^{\infty} 
L_x^3} \|u_5\|_{L_t^2 L_x^6} \|u_3\|_{L_{xt}^{\infty}} 
\|u_6\|_{L_{xt}^{\infty}} 
\\
& \lesssim & \frac{1}{\prod_{i=1}^6 N_i^{\frac{3}{4}} N^{\frac{3}{2}}} 
N_2^{\frac{1}{2}} N_4^{\frac{1}{2}} N_3^{\frac{3}{2}+} N_6^{\frac{3}{2}+} 
\prod_{i=1}^6 \|u_i\|_{X^{0,\frac{1}{2}+}} \\ 
& \lesssim & \frac{N_3^{\frac{3}{4}} N_6^{\frac{3}{4}+}}{N_1^{\frac{3}{4}} 
N_2^{\frac{1}{4}} N_4^{\frac{1}{4}} N_5^{\frac{3}{4}} N^{\frac{3}{2}}} 
\prod_{i=1}^6 \|u_i\|_{X^{0,\frac{1}{2}+}}\\
& \lesssim & \frac{1 \wedge N_{min}^{0+}}{N_1^{0+} N_2^{\frac{1}{4}-} N_4^{0+} 
N_5^{\frac{1}{4}-} N^{\frac{3}{2}}} \prod_{i=1}^6 
\|u_i\|_{X^{0,\frac{1}{2}+}}\\
& \lesssim & \frac{1 \wedge N_{min}^{0+}}{ N_{max}^{0+} N^{2-}} \prod_{i=1}^6 
\|u_i\|_{X^{0,\frac{1}{2}+}} \, .
\end{eqnarray*}
{\bf b.} $N_1 \ge N_2 \ge N \ge N_3$.\\
This can be handled similarly as case a. without the factor 
$(\frac{N_3}{N})^{\frac{1}{4}}$. Thus
\begin{eqnarray*}
C & \lesssim & \frac{N_3^{\frac{1}{2}+} N_6^{\frac{3}{4}+}}{N_1^{\frac{3}{4}} 
N_2^{\frac{1}{4}} N_4^{\frac{1}{4}} N_5^{\frac{3}{4}} N^{\frac{5}{4}}} 
\prod_{i=1}^6 \|u_i\|_{X^{0,\frac{1}{2}+}}
 \lesssim  \frac{1 \wedge N_{min}^{0+}}{N_1^{0+} N_2^{\frac{1}{2}-} N_4^{0+} 
N_5^{\frac{1}{4}-} N^{\frac{5}{4}}} \prod_{i=1}^6 
\|u_i\|_{X^{0,\frac{1}{2}+}}\\
& \lesssim & \frac{1 \wedge N_{min}^{0+}}{ N_{max}^{0+} N^{2-}} \prod_{i=1}^6 
\|u_i\|_{X^{0,\frac{1}{2}+}} \, .
\end{eqnarray*}
{\bf c.} $N_1 \ge N \ge N_2 \ge N_3$.\\
This can be handled as case b. without the factor 
$(\frac{N_2}{N})^{\frac{1}{4}}$ leading to the bound
\begin{eqnarray*}
C & \lesssim & \frac{N_3^{\frac{1}{2}+} N_6^{\frac{3}{4}+}}{N_1^{\frac{3}{4}} 
N_2^{\frac{1}{2}} N_4^{\frac{1}{4}} N_5^{\frac{3}{4}} N} \prod_{i=1}^6 
\|u_i\|_{X^{0,\frac{1}{2}+}}
 \lesssim  \frac{1 \wedge N_{min}^{0+}}{N_1^{\frac{3}{4}-} N_4^{\frac{1}{4}-} 
N} 
\prod_{i=1}^6 \|u_i\|_{X^{0,\frac{1}{2}+}}\\
& \lesssim & \frac{1 \wedge N_{min}^{0+}}{ N_{max}^{0+} N^{2-}} \prod_{i=1}^6 
\|u_i\|_{X^{0,\frac{1}{2}+}} \, .
\end{eqnarray*} 
{\bf Case 2:} $N_4 \ge N_5 \ge N \ge N_6$.\\
{\bf a.} $N_1 \ge N_2 \ge N_3 \ge N$. \\
This is handled like case 1a. without the factor  
$(\frac{N_6}{N})^{\frac{1}{4}}$ and gives
\begin{eqnarray*}
C & \lesssim & \frac{N_3^{\frac{3}{4}+} N_6^{\frac{1}{2}+}}{N_1^{\frac{3}{4}} 
N_2^{\frac{1}{4}} N_4^{\frac{1}{4}} N_5^{\frac{3}{4}} N^{\frac{5}{4}}} 
\prod_{i=1}^6 \|u_i\|_{X^{0,\frac{1}{2}+}}
 \lesssim  \frac{1 \wedge N_{min}^{0+}}{N_1^{0+} N_2^{\frac{1}{4}-} 
N_4^{0+}N_5^{\frac{1}{2}-} N^{\frac{5}{4}}} \prod_{i=1}^6 
\|u_i\|_{X^{0,\frac{1}{2}+}}\\
& \lesssim & \frac{1 \wedge N_{min}^{0+}}{ N_{max}^{0+} N^{2-}} \prod_{i=1}^6 
\|u_i\|_{X^{0,\frac{1}{2}+}} \, .
\end{eqnarray*} 
{\bf b.} $N_1 \ge N_2 \ge N \ge N_3$.\\
It can be handled like case 1b. without the factor   
$(\frac{N_6}{N})^{\frac{1}{4}}$ leading to
\begin{eqnarray*}
C & \lesssim & \frac{N_3^{\frac{1}{2}+} N_6^{\frac{1}{2}+}}{N_1^{\frac{3}{4}} 
N_2^{\frac{1}{4}} N_4^{\frac{1}{4}} N_5^{\frac{3}{4}} N} \prod_{i=1}^6 
\|u_i\|_{X^{0,\frac{1}{2}+}}
 \lesssim  \frac{1 \wedge N_{min}^{0+}}{N_1^{0+} N_2^{\frac{1}{2}-} 
N_4^{0+}N_5^{\frac{1}{2}-} N} \prod_{i=1}^6 \|u_i\|_{X^{0,\frac{1}{2}+}}\\
& \lesssim & \frac{1 \wedge N_{min}^{0+}}{ N_{max}^{0+} N^{2-}} \prod_{i=1}^6 
\|u_i\|_{X^{0,\frac{1}{2}+}} \, .
\end{eqnarray*}
{\bf c.} $N_1 \ge N \ge N_2 \ge N_3$.\\
As in case 1c. without the factor   $(\frac{N_6}{N})^{\frac{1}{4}}$ we get
\begin{eqnarray*}
C & \lesssim & \frac{N_3^{\frac{1}{2}+} N_6^{\frac{1}{2}+}}{N_1^{\frac{3}{4}} 
N_2^{\frac{1}{2}} N_4^{\frac{1}{4}} N_5^{\frac{3}{4}} N^{\frac{3}{4}}} 
\prod_{i=1}^6 \|u_i\|_{X^{0,\frac{1}{2}+}}
 \lesssim  \frac{1 \wedge N_{min}^{0+}}{N_1^{\frac{3}{4}-} 
N_4^{\frac{1}{4}}N_5^{\frac{1}{4}-} N^{\frac{3}{4}}} \prod_{i=1}^6 
\|u_i\|_{X^{0,\frac{1}{2}+}}\\
& \lesssim & \frac{1 \wedge N_{min}^{0+}}{ N_{max}^{0+} N^{2-}} \prod_{i=1}^6 
\|u_i\|_{X^{0,\frac{1}{2}+}} \, .
\end{eqnarray*}
{\bf Case 3:} $N_4 \gtrsim N \ge N_5 \ge N_6$.\\
{\bf a.} $N_1 \ge N_2 \ge N_3 \ge N$.\\
We have $$ M(\xi_1,...,\xi_6) \lesssim \prod_{i=1}^3 
(\frac{N_i}{N})^{\frac{1}{4}} \prod_{i=1}^6 N_i^{-1} \, . $$
Thus by Strichartz and Sobolev and the bilinear Strichartz estimate (\ref{B2}) 
we get
\begin{eqnarray*}
C & \lesssim &  \prod_{i=1}^3 (\frac{N_i}{N})^{\frac{1}{4}} \prod_{i=1}^6 
N_i^{-1} \|u_1\|_{L_t^2 L_x^6} \|u_2\|_{L_t^{\infty} L_x^3}\|u_3 u_4\|_{L_{tx}^2} 
\|u_5\|_{L_{xt}^{\infty}} \|u_6\|_{L_{xt}^{\infty}} \\
& \lesssim &  \prod_{i=1}^3 (\frac{N_i}{N})^{\frac{1}{4}} \prod_{i=1}^6 
N_i^{-1} 
\frac{N_2^{\frac{1}{2}} N_3^{1+}}{N_4^{\frac{1}{2}-}} (N_5^{\frac{3}{2}+} + 
N_5^{\frac{3}{2}-})  (N_6^{\frac{3}{2}+} + N_6^{\frac{3}{2}-})\prod_{i=1}^6 
\|u_i\|_{X^{0,\frac{1}{2}+}} \\
& \lesssim & \frac{N_3^{\frac{1}{4}+} (N_5^{\frac{1}{2}+} + 
N_5^{\frac{1}{2}-})  
(N_6^{\frac{1}{2}+} + N_6^{\frac{1}{2}-})}{N_1^{\frac{3}{4}} N_2^{\frac{1}{4}} 
N_4^{\frac{3}{2}-} N^{\frac{3}{4}}} \prod_{i=1}^6 \|u_i\|_{X^{0,\frac{1}{2}+}} 
\\
& \lesssim & \frac{N^{1+} (1 \wedge N_{min}^{0+})}{ N_{max}^{0+} N^{3-}} 
\prod_{i=1}^6 \|u_i\|_{X^{0,\frac{1}{2}+}} \\
& \lesssim & \frac{1 \wedge N_{min}^{0+}}{ N_{max}^{0+} N^{2-}} \prod_{i=1}^6 
\|u_i\|_{X^{0,\frac{1}{2}+}} \, .
\end{eqnarray*}
{\bf b.} $ N_1 \ge N_2 \ge N \ge N_3$.\\
Similarly as in case a. without the factor $(\frac{N_3}{N})^{\frac{1}{4}}$ we 
get
\begin{eqnarray*}
C & \lesssim & \frac{N_3^{0+} (N_5^{\frac{1}{2}+} + N_5^{\frac{1}{2}-})  
(N_6^{\frac{1}{2}+} + N_6^{\frac{1}{2}-})}{N_1^{\frac{3}{4}} N_2^{\frac{1}{4}} 
N_4^{\frac{3}{2}-} N^{\frac{1}{2}}} \prod_{i=1}^6 \|u_i\|_{X^{0,\frac{1}{2}+}} 
\\
& \lesssim & \frac{N^{1+} (1 \wedge N_{min}^{0+})}{ N_{max}^{0+} N^{3-}} 
\prod_{i=1}^6 \|u_i\|_{X^{0,\frac{1}{2}+}} \\
& \lesssim & \frac{1 \wedge N_{min}^{0+}}{ N_{max}^{0+} N^{2-}} \prod_{i=1}^6 
\|u_i\|_{X^{0,\frac{1}{2}+}} \, .
\end{eqnarray*}
{\bf c.} $N_1 \ge N \ge N_2 \ge N_3$.\\
The multiplier is estimated as follows using the mean value theorem
\begin{eqnarray*}
M(\xi_1,...,\xi_6) & \lesssim & 
|\frac{m(\xi_1+\xi_2+\xi_3)-m(\xi_1)}{m(\xi_1)}| 
\prod_{i=1}^6 |\xi_i|^{-1} \\ 
& \sim & \frac{|(\nabla m)(\xi_1)(\xi_2 + \xi_3)|}{|m(\xi_1)|} \prod_{i=1}^6 
|\xi_i|^{-1} 
 \lesssim  \frac{N_2}{N_1} \prod_{i=1}^6 |\xi_i|^{-1} \, . 
\end{eqnarray*} 
Thus we get by Sobolev and Strichartz:
\begin{eqnarray*}
C & \lesssim & \frac{N_2}{N_1} \prod_{i=1}^6 N_i^{-1} \|u_1\|_{L_t^2 L_x^6} 
\|u_2\|_{L_t^{\infty} L_x^6}\|u_3\|_{L_t^{\infty} L_x^{6+}} \|u_4\|_{L_t^2 
L_x^{6-}} \|u_5\|_{L_t^{\infty} L_x^6} \|u_6\|_{L_t^{\infty} L_x^{6+}} \\
& \lesssim & \frac{N_2}{ N_1}  \prod_{i=1}^6 N_i^{-1} N_2 N_3^{1+} N_5 
N_6^{1+} 
\prod_{i=1}^6 \|u_i\|_{X^{0,\frac{1}{2}+}} \\
& \lesssim & \frac{N_2 N_3^{0+} N_6^{0+}}{N_1^2 N_4} \prod_{i=1}^6 
\|u_i\|_{X^{0,\frac{1}{2}+}} \\
& \lesssim & \frac{1 \wedge N_{min}^{0+}}{ N_{max}^{0+} N^{2-}} \prod_{i=1}^6 
\|u_i\|_{X^{0,\frac{1}{2}+}} \, .
\end{eqnarray*}
{\bf Case 4:} $N \gg N_4 \ge N_5 \ge N_6$.\\
{\bf a.} $N_1 \ge N_2 \ge N_3 \ge N$.\\
We have
\begin{eqnarray*}
C & \lesssim & (\frac{N_1}{N})^{\frac{1}{4}} (\frac{N_2}{N})^{\frac{1}{4}} 
(\frac{N_3}{N})^{\frac{1}{4}}  \prod_{i=1}^6 N_i^{-1} \|u_1\|_{L_t^2 L_x^{6-}} 
\|u_2\|_{L_t^2 L_x^{6-}} \prod_{i=3}^6 \|u_i\|_{L_t^{\infty} L_x^{6+}} \\
& \lesssim &  (\frac{N_1}{N})^{\frac{1}{4}} (\frac{N_2}{N})^{\frac{1}{4}} 
(\frac{N_3}{N})^{\frac{1}{4}}  \prod_{i=1}^6 N_i^{-1} N_3^{1+} N_4^{1+} 
N_5^{1+} 
N_6^{1+} \prod_{i=1}^6 \|u_i\|_{X^{0,\frac{1}{2}+}} \\
& \lesssim & \frac{ N_3^{\frac{1}{4}+} N_4^{0+} N_5^{0+} 
N_6^{0+}}{N_1^{\frac{3}{4}} N_2^{\frac{3}{4}} N^{\frac{3}{4}}} \prod_{i=1}^6 
\|u_i\|_{X^{0,\frac{1}{2}+}} \\
& \lesssim & \frac{1 \wedge N_{min}^{0+}}{ N_{max}^{0+} N^{2-}} \prod_{i=1}^6 
\|u_i\|_{X^{0,\frac{1}{2}+}} \, .
\end{eqnarray*}
{\bf b.} $N_1 \ge N_2 \gtrsim N \gtrsim N_3$.\\
We get the same bound as in case a. without the factor 
$(\frac{N_3}{N})^{\frac{1}{4}}$ leading to the same estimate. \\
{\bf c.} $N_1 \ge N \gg N_2 \ge N_3$.\\
This case cannot occur, because $\sum_{i=1}^6 \xi_i = 0$.\\
{\bf 4.} Next we aim to show
$$ \int_0^{\delta}| \langle I(u^3) - (Iu)^3, I(u^2) \rangle| dt \lesssim 
N^{-\frac{5}{2}+} \|\nabla Iu\|_{X^{0,\frac{1}{2}+}}^5 \, , $$
which follows from
\begin{eqnarray}
\label{15}
D  =  |\int_0^{\delta} \int_* M(\xi_1,...,\xi_5) \prod_{i=1}^5 
\widehat{u_i}(\xi_i,t)d\xi_1 ...d\xi_5 dt| \lesssim N^{-\frac{5}{2}+} 
\prod_{i=1}^5 \|u_i\|_{X^{0,\frac{1}{2}+}} \, , 
\end{eqnarray}
where
$$ M(\xi_1,...,\xi_5) := 
\frac{|m(\xi_1+\xi_2+\xi_3)-m(\xi_1)m(\xi_2)m(\xi_3)|}{m(\xi_1)m(\xi_2)m(\xi_3)} 
\frac{m(\xi_4+\xi_5)}{m(\xi_4)m(\xi_5)} \prod_{i=1}^5 |\xi_i|^{-1} \, . $$
Assume w.l.o.g. $N_1\ge  N_2 \ge N_3$ , $N_1 \ge N$ and $N_4 \ge N_5$.\\
{\bf Case 1:} $N_4 \ge N_5 \ge N$.\\
{\bf a.} $N_1 \ge N_2 \ge N_3 \ge N$.\\
Thus by Strichartz and Sobolev:
\begin{eqnarray*}
D & \lesssim &  \prod_{i=1}^5 (\frac{N_i}{N})^{\frac{1}{4}} \prod_{i=1}^5 
N_i^{-1} \|u_1\|_{L_t^{\infty} L_x^3} \|u_2\|_{L_t^2 L_x^6} 
\|u_3\|_{L^{\infty} 
L_x^{6+}}\|u_4\|_{L_t^2 L_x^{6-}} \|u_5\|_{L_t^{\infty} L_x^{6+}} \\
& \lesssim &  \prod_{i=1}^5 (\frac{N_i}{N})^{\frac{1}{4}} \prod_{i=1}^5 
N_i^{-1} 
N_1^{\frac{1}{2}} N_3^{1+} N_5^{1+} \prod_{i=1}^5 \|u_i\|_{X^{0,\frac{1}{2}+}} 
\\ 
& \lesssim & \frac{N_3^{\frac{1}{4}+} N_5^{\frac{1}{4}+}}{N_1^{\frac{1}{4}} 
N_2^{\frac{3}{4}} N_4^{\frac{3}{4}} N^{\frac{5}{4}}} \prod_{i=1}^5 
\|u_i\|_{X^{0,\frac{1}{2}+}}\\
& \lesssim & \frac{1 \wedge N_{min}^{0+}}{ N_{max}^{0+} N^{\frac{5}{2}-}} 
\prod_{i=1}^5 \|u_i\|_{X^{0,\frac{1}{2}+}} \, .
\end{eqnarray*}
{\bf b.} $N_1 \ge N_2 \ge N \ge N_3$.\\
As in case a. (without the factor $(\frac{N_3}{N})^{\frac{1}{4}}$) we get the 
same estimate. \\
{\bf c.} $N_1 \ge N \ge N_2 \ge N_3$ and $N_1 \gg N_2$.\\
The multiplier is estimated as follows:
\begin{eqnarray*}
|\frac{m(\xi_1+\xi_2+\xi_3)-m(\xi_1) m(\xi_2) m(\xi_3)}{m(\xi_1) m(\xi_2) 
m(\xi_3)}|&  = & |\frac{m(\xi_1+\xi_2+\xi_3)-m(\xi_1)}{m(\xi_1)}|\\ &   \sim & 
|\frac{(\nabla m)(\xi_1) (\xi_2 + \xi_3)}{m(\xi_1)}| \lesssim \frac{N_2}{N_1} 
\, 
. 
\end{eqnarray*}
Thus as in case a. we get
\begin{eqnarray*}
D & \lesssim & \frac{N_2}{N_1} (\frac{N_4}{N})^{\frac{1}{4}} 
(\frac{N_5}{N})^{\frac{1}{4}} \prod_{i=1}^5 N_i^{-1} N_1^{\frac{1}{2}} 
N_3^{1+} 
N_5^{1+} \prod_{i=1}^5 \|u_i\|_{X^{0,\frac{1}{2}+}} \\ 
& \lesssim & \frac{N_3^{0+} N_5^{\frac{1}{4}+}}{N_1^{\frac{3}{2}} 
N_4^{\frac{3}{4}} N^{\frac{1}{2}}} \prod_{i=1}^5 
\|u_i\|_{X^{0,\frac{1}{2}+}}\\
& \lesssim & \frac{1 \wedge N_{min}^{0+}}{ N_{max}^{0+} N^{\frac{5}{2}-}} 
\prod_{i=1}^5 \|u_i\|_{X^{0,\frac{1}{2}+}} \, .
\end{eqnarray*}
{\bf Case 2:} $N_4 \gtrsim N \ge N_5$.\\
{\bf a.} $N_1 \ge N_2 \ge N_3 \ge N$\\
As in case 1a. (without the factor  $(\frac{N_5}{N})^{\frac{1}{4}}$) we get 
the 
same estimate.\\
{\bf b.} $N_1 \ge N_2 \ge N \ge N_3$.\\
As in case 1a.  (without the factor  
$(\frac{N_3}{N})^{\frac{1}{4}}(\frac{N_5}{N})^{\frac{1}{4}}$) we get the same 
estimate.\\
{\bf c.} $N_1 \ge N \ge N_2 \ge N_3$.\\
As in case 1c.  (without the factor  $(\frac{N_5}{N})^{\frac{1}{4}}$) we get 
the 
same estimate.\\
{\bf Case 3:} $ N \gg N_4 \ge N_5$ $\Rightarrow N_1 \sim N_2 \gtrsim N$, 
because 
$\sum_{i=1}^5 \xi_i = 0$. \\
{\bf a.} $N_1 \sim N_2 \ge N_3 \ge N$. \\
By Strichartz and Sobolev and the bilinear Strichartz estimate (\ref{B2}) we 
get
\begin{eqnarray*}
D & \lesssim &  \prod_{i=1}^3 (\frac{N_i}{N})^{\frac{1}{4}} \prod_{i=1}^5 
N_i^{-1} \|u_1 u_3\|_{L_{xt}^2} \|u_2\|_{L_t^2 L_x^{6-}} \|u_4\|_{L_t^{\infty} 
L_x^{6}} \|u_5\|_{L_t^{\infty} L_x^{6+}} \\
& \lesssim &  \prod_{i=1}^3 (\frac{N_i}{N})^{\frac{1}{4}} \prod_{i=1}^5 
N_i^{-1} 
\frac{ N_3^{1+}}{ N_1^{\frac{1}{2}-}}N_4 N_5^{1+} \prod_{i=1}^5 
\|u_i\|_{X^{0,\frac{1}{2}+}} \\
& \lesssim & \frac{N_3^{\frac{1}{4}+} N_5^{0+}}{N_1^{\frac{5}{4}-} 
N_2^{\frac{3}{4}} N^{\frac{3}{4}}} \prod_{i=1}^5 
\|u_i\|_{X^{0,\frac{1}{2}+}}
 \lesssim  \frac{1 \wedge N_{min}^{0+}}{ N_{max}^{0+} N^{\frac{5}{2}-}} 
\prod_{i=1}^5 \|u_i\|_{X^{0,\frac{1}{2}+}} \, .
\end{eqnarray*}
{\bf b.} $N_1 \sim N_2 \gtrsim N \ge N_3$.\\
Without the factor  $(\frac{N_3}{N})^{\frac{1}{4}}$) we get the same estimate 
as 
in case a. \\
{\bf 5.} Next we want to show
$$ \int_0^{\delta}| \langle I(u^2)-(Iu)^2, I(u^3) \rangle| dt \lesssim 
N^{-\frac{5}{2}+} \|\nabla Iu\|_{X^{0,\frac{1}{2}+}}^5 \, . $$
We have to prove
\begin{eqnarray}
\label{16}
E  =  |\int_0^{\delta} \int_* M(\xi_1,...,\xi_5) \prod_{i=1}^5 
\widehat{u_i}(\xi_i,t)d\xi_1 ...d\xi_5 dt| \lesssim N^{-\frac{5}{2}+} 
\prod_{i=1}^5 \|u_i\|_{X^{0,\frac{1}{2}+}} \, , 
\end{eqnarray}
where
$$ M(\xi_1,...,\xi_5) := 
\frac{|m(\xi_1+\xi_2)-m(\xi_1)m(\xi_2)|}{m(\xi_1)m(\xi_2)} 
\frac{m(\xi_3+\xi_4+\xi_5)}{m(\xi_3)m(\xi_4)m(\xi_5)} \prod_{i=1}^5 
|\xi_i|^{-1} 
\, . $$
Assume w.l.o.g. $N_3\ge  N_4 \ge N_5$ , $N_1 \ge N_2$ and $N_1 \ge N$.\\
{\bf Case 1:} $N_3 \ge N_4 \ge N_5 \ge N$. \\
{\bf a.} $N_1 \ge N_2 \gtrsim N$.\\
This case can be handled exactly as in 4. case 1a.\\
{\bf b.} $N_1 \ge N \gg N_2$.\\
We get by use of the mean value theorem for the first fraction and estimating 
similarly as in 4. case 1a. (interchanging the roles of $u_2$ and $u_3$):
\begin{eqnarray*}
E & \lesssim & \frac{N_2}{N_1} \prod_{i=3}^5 (\frac{N_i}{N})^{\frac{1}{4}} 
\prod_{i=1}^5 N_i^{-1} N_1^{\frac{1}{2}} N_2^{1+} N_5^{1+} \prod_{i=1}^5 
\|u_i\|_{X^{0,\frac{1}{2}+}} \\ 
& \lesssim & \frac{N_2^{1+} N_5^{\frac{1}{4}+}}{N_1^{\frac{3}{2}} 
N_3^{\frac{3}{4}} N_4^{\frac{3}{4}} N^{\frac{3}{4}}} \prod_{i=1}^5 
\|u_i\|_{X^{0,\frac{1}{2}+}}\\
& \lesssim & \frac{1 \wedge N_{min}^{0+}}{ N_{max}^{0+} N^{\frac{5}{2}-}} 
\prod_{i=1}^5 \|u_i\|_{X^{0,\frac{1}{2}+}} \, .
\end{eqnarray*}
{\bf Case 2:} $N_3 \ge N_4 \ge N \ge N_5$.\\
{\bf a.} $N_1 \ge N_2 \gtrsim N$.\\
By Strichartz and Sobolev we get
\begin{eqnarray*}
E & \lesssim &  \prod_{i=1}^4 (\frac{N_i}{N})^{\frac{1}{4}} \prod_{i=1}^5 
N_i^{-1} \|u_1\|_{L_t^{\infty} L_x^3}\|u_2\|_{L_t^2 L_x^{6-}} \|u_3\|_{L_t^2 
L_x^{6}} \|u_4\|_{L_t^{\infty} L_x^{6}} \|u_5\|_{L_t^{\infty} L_x^{6+}} \\
& \lesssim &  \prod_{i=1}^4 (\frac{N_i}{N})^{\frac{1}{4}} \prod_{i=1}^5 
N_i^{-1} 
N_1^{\frac{1}{2}} N_4 N_5^{1+} \prod_{i=1}^5 \|u_i\|_{X^{0,\frac{1}{2}+}} \\ 
& \lesssim & \frac{N_4^{\frac{1}{4}} N_5^{0+}}{N_1^{\frac{1}{4}} 
N_2^{\frac{3}{4}}N_3^{\frac{3}{4}} N} \prod_{i=1}^5 
\|u_i\|_{X^{0,\frac{1}{2}+}}\\
& \lesssim & \frac{1 \wedge N_{min}^{0+}}{ N_{max}^{0+} N^{\frac{5}{2}-}} 
\prod_{i=1}^5 \|u_i\|_{X^{0,\frac{1}{2}+}} \, .
\end{eqnarray*}
{\bf b.} $N_1 \ge N_2 \gg N_2$.\\
By the mean value theorem we get as in a. (slightly modified):
\begin{eqnarray*}
E & \lesssim & \frac{N_2}{N_1} (\frac{N_3}{N})^{\frac{1}{4}}  
(\frac{N_4}{N})^{\frac{1}{4}} \prod_{i=1}^5 N_i^{-1} N_1^{\frac{1}{2}}N_2^{0+} 
N_4^{1-} N_5^{1+} \prod_{i=1}^5 \|u_i\|_{X^{0,\frac{1}{2}+}} \\ 
& \lesssim & \frac{N_2^{0+} N_4^{\frac{1}{4}-} N_5^{0+}}{N_1^{\frac{3}{2}} 
N_3^{\frac{3}{4}} N^{\frac{1}{2}}} \prod_{i=1}^5 
\|u_i\|_{X^{0,\frac{1}{2}+}}\\
& \lesssim & \frac{1 \wedge N_{min}^{0+}}{ N_{max}^{0+} N^{\frac{5}{2}-}} 
\prod_{i=1}^5 \|u_i\|_{X^{0,\frac{1}{2}+}} \, .
\end{eqnarray*}
{\bf Case 3:} $N_3 \gtrsim N \ge N_4 \ge N_5$.\\
{\bf a.} $N_1 \ge N_2 \ge N$.\\
The second fraction is bounded, so that as in case 2a. we get
\begin{eqnarray*}
E & \lesssim &  (\frac{N_1}{N})^{\frac{1}{4}} (\frac{N_2}{N})^{\frac{1}{4}} 
\prod_{i=1}^5 N_i^{-1}  N_1^{\frac{1}{2}} N_4 N_5^{1+} \prod_{i=1}^5 
\|u_i\|_{X^{0,\frac{1}{2}+}} \\
& \lesssim & \frac{ N_5^{0+}}{N_1^{\frac{1}{4}} N_2^{\frac{3}{4}} N_3 
N^{\frac{1}{2}}} \prod_{i=1}^5 \|u_i\|_{X^{0,\frac{1}{2}+}}\\
& \lesssim & \frac{1 \wedge N_{min}^{0+}}{ N_{max}^{0+} N^{\frac{5}{2}-}} 
\prod_{i=1}^5 \|u_i\|_{X^{0,\frac{1}{2}+}} \, . 
\end{eqnarray*}
{\bf b.} $N_1 \ge N \ge N_2$.\\
Using the mean value theorem and interchanging the roles of $u_1$ and $u_2$ in 
case a. gives 
\begin{eqnarray*}
E & \lesssim & \frac{N_2}{N_1} \prod_{i=1}^5 N_i^{-1} 
N_2^{\frac{1}{2}} N_4 N_5^{1+} \prod_{i=1}^5 \|u_i\|_{X^{0,\frac{1}{2}+}} \\
 & \lesssim & \frac{N_2^{\frac{1}{2}} N_5^{0+}}{N_1^2 N_3} \prod_{i=1}^5 
\|u_i\|_{X^{0,\frac{1}{2}+}}\\
& \lesssim & \frac{1 \wedge N_{min}^{0+}}{ N_{max}^{0+} N^{\frac{5}{2}-}} 
\prod_{i=1}^5 \|u_i\|_{X^{0,\frac{1}{2}+}} \, .
\end{eqnarray*}
{\bf Case 4:} $N \gg N_3 \ge N_4 \ge N_5$ $\Rightarrow N_1 \sim N_2 \gtrsim 
N$, 
because $\sum_{i=1}^5 \xi_i = 0$. \\
This gives
by Strichartz, Sobolev and the bilinear Strichartz estimate (\ref{B1}) we get
\begin{eqnarray*}
E & \lesssim & (\frac{N_1}{N})^{\frac{1}{4}} (\frac{N_2}{N})^{\frac{1}{4}} 
\prod_{i=1}^5 N_i^{-1} \|u_1 u_3\|_{L_{xt}^2} \|u_2\|_{L_t^2 L_x^{6-}} 
\|u_4\|_{L_t^{\infty} L_x^{6}} \|u_5\|_{L_t^{\infty} L_x^{6+}} \\
& \lesssim &  (\frac{N_1}{N})^{\frac{1}{4}} (\frac{N_2}{N})^{\frac{1}{4}} 
\prod_{i=1}^5 N_i^{-1} \frac{N_3}{N_1^{\frac{1}{2}}} N_4 N_5^{1+} 
\prod_{i=1}^5 
\|u_i\|_{X^{0,\frac{1}{2}+}} \\ 
& \lesssim & \frac{ N_5^{0+}}{N_1^{\frac{5}{4}} N_2^{\frac{3}{4}} 
N^{\frac{1}{2}}} \prod_{i=1}^5 \|u_i\|_{X^{0,\frac{1}{2}+}}\\
& \lesssim & \frac{1 \wedge N_{min}^{0+}}{ N_{max}^{0+} N^{\frac{5}{2}-}} 
\prod_{i=1}^5 \|u_i\|_{X^{0,\frac{1}{2}+}} \, .
\end{eqnarray*}
{\bf 6.} Next we want to prove
$$ \int_0^{\delta} |\langle I(u^3) - (Iu)^3, Iu \rangle | dt \lesssim N^{-3+} 
\|\nabla Iu\|_{X^{0,\frac{1}{2}+}}^4 \, , $$
which follows from
\begin{eqnarray}
\label{17}
F  =  |\int_0^{\delta} \int_* M(\xi_1,...,\xi_4) \prod_{i=1}^4 
\widehat{u_i}(\xi_i,t)d\xi_1 ...d\xi_4 dt| \lesssim N^{-3+} \prod_{i=1}^4 
\|u_i\|_{X^{0,\frac{1}{2}+}} \, , 
\end{eqnarray}
where
$$ M(\xi_1,\xi_2,\xi_3) := 
\frac{|m(\xi_1+\xi_2+\xi_3)-m(\xi_1)m(\xi_2)m(\xi_3)|}{m(\xi_1)m(\xi_2)m(\xi_3)} 
\prod_{i=1}^4 |\xi_i|^{-1} \, . $$
Assume w.l.o.g. $N_1\ge  N_2 \ge N_3$ and $N_1 \ge N$.\\
{\bf Case 1:} $ N \gg N_4 $.\\
{\bf a.} $N_1 \ge N_2 \ge N_3 \ge N$.\\
By the bilinear Strichartz refinement (\ref{B2}) we get
\begin{eqnarray*}
E & \lesssim & \prod_{i=1}^3(\frac{N_i}{N})^{\frac{1}{4}} \prod_{i=1}^4 
N_i^{-1} 
\|u_2 u_4\|_{L_{xt}^2}\|u_1 u_3\|_{L_{xt}^2}  \\
& \lesssim & \prod_{i=1}^3(\frac{N_i}{N})^{\frac{1}{4}} \prod_{i=1}^4 N_i^{-1} 
\frac{N_4^{1+}}{N_2^{\frac{1}{2}-}} \frac{N_3^{1+}}{N_1^{\frac{1}{2}-}}        
\prod_{i=1}^4 \|u_i\|_{X^{0,\frac{1}{2}+}} \\ 
& \lesssim & \frac{N_3^{\frac{1}{4}+} N_4^{0+}}{N_1^{\frac{5}{4}-} 
N_2^{\frac{5}{4}-} N^{\frac{3}{4}}} \prod_{i=1}^4 
\|u_i\|_{X^{0,\frac{1}{2}+}}
 \lesssim  \frac{1 \wedge N_{min}^{0+}}{ N_{max}^{0+} N^{3-}} \prod_{i=1}^4 
\|u_i\|_{X^{0,\frac{1}{2}+}} \, .
\end{eqnarray*}
{\bf b.} $N_1 \ge N_2 \gtrsim N \ge N_3$.\\
Similarly as in case a. we get 
\begin{eqnarray*}
E & \lesssim & \prod_{i=1}^2(\frac{N_i}{N})^{\frac{1}{4}} \prod_{i=1}^4 
N_i^{-1} 
\frac{N_4^{1+}}{N_2^{\frac{1}{2}-}} \frac{N_3^{1+}}{N_1^{\frac{1}{2}-}}        
\prod_{i=1}^4 \|u_i\|_{X^{0,\frac{1}{2}+}} \\ 
& \lesssim & \frac{N_3^{0+} N_4^{0+}}{N_1^{\frac{5}{4}-} N_2^{\frac{5}{4}-} 
N^{\frac{1}{2}}} \prod_{i=1}^4 \|u_i\|_{X^{0,\frac{1}{2}+}}
 \lesssim \frac{1 \wedge N_{min}^{0+}}{ N_{max}^{0+} N^{3-}} \prod_{i=1}^4 
\|u_i\|_{X^{0,\frac{1}{2}+}} \, .
\end{eqnarray*}
{\bf c.} $N_1 \ge N \gg N_2 \ge N_3$.\\
This case does not occur, because $\sum_{i=1}^4 \xi_i = 0$.\\
{\bf Case 2:} $N_4 \gtrsim N$.\\
{\bf a.} $N_1 \ge N_2 \ge N_3 \ge N$ and\\
{\bf b.} $N_1 \ge N_2 \gtrsim N \ge N_3$ can be treated as in case 1.\\
{\bf c.} $N_1 \ge N \gg N_2 \ge N_3$ $\Rightarrow N_1 \sim N_4$.\\
By the mean value theorem and the bilinear Strichartz refinement  (\ref{B2}) 
we 
get
\begin{eqnarray*}
E & \lesssim & \frac{N_2}{N_1} \prod_{i=1}^4 N_i^{-1} \|u_1 
u_2\|_{L_{xt}^2}\|u_3 u_4\|_{L_{xt}^2}  \\
& \lesssim & \frac{N_2}{N_1} \prod_{i=1}^4 N_i^{-1} 
\frac{N_2}{N_1^{\frac{1}{2}}} \frac{N_3^{1+}}{N_4^{\frac{1}{2}-}}        
\prod_{i=1}^4 \|u_i\|_{X^{0,\frac{1}{2}+}} \\
& \lesssim & \frac{N_2 N_3^{0+}}{N_1^{\frac{5}{2}} N_4^{\frac{3}{2}-} }
\prod_{i=1}^4 \|u_i\|_{X^{0,\frac{1}{2}+}}
 \lesssim  \frac{1 \wedge N_{min}^{0+}}{ N_{max}^{0+} N^{3-}} \prod_{i=1}^4 
\|u_i\|_{X^{0,\frac{1}{2}+}} \, .
\end{eqnarray*}
{\bf 7.} Next we prove
$$ \int_0^{\delta} |\langle I(u^2) - (Iu)^2, (Iu)^2 \rangle | dt \lesssim 
N^{-3+} \|\nabla Iu\|_{X^{0,\frac{1}{2}+}}^4 \, , $$
which follows from
\begin{eqnarray}
\label{18}
F  =  |\int_0^{\delta} \int_* M(\xi_1,...,\xi_4) \prod_{i=1}^4 
\widehat{u_i}(\xi_i,t)d\xi_1 ...d\xi_4 dt| \lesssim N^{-3+} \prod_{i=1}^4 
\|u_i\|_{X^{0,\frac{1}{2}+}} \, , 
\end{eqnarray}
where
$$ M(\xi_1,\xi_2,\xi_3) := 
\frac{|m(\xi_1+\xi_2)-m(\xi_1)m(\xi_2)|}{m(\xi_1)m(\xi_2)}  
\frac{m(\xi_3+\xi_4)}{m(\xi_3)m(\xi_4)} \prod_{i=1}^4 |\xi_i|^{-1} \, . $$
Assume w.l.o.g. $N_1\ge  N_2$, $N_1 \ge N$ and $N_3 \ge N_4$.\\
{\bf Case 1:} $N_3 \ge N_4 \gtrsim N$.\\
{\bf a.} $N_1 \ge N_2 \ge N$.\\
By the bilinear Strichartz refinement (\ref{B2}) we get
\begin{eqnarray*}
F & \lesssim & \prod_{i=1}^4(\frac{N_i}{N})^{\frac{1}{4}} \prod_{i=1}^4 
N_i^{-1} 
\|u_1 u_2\|_{L_{xt}^2}\|u_3 u_4\|_{L_{xt}^2}  \\
& \lesssim & \prod_{i=1}^4(\frac{N_i}{N})^{\frac{1}{4}} \prod_{i=1}^4 N_i^{-1} 
\frac{N_2^{1+}}{N_1^{\frac{1}{2}-}} \frac{N_4^{1+}}{N_3^{\frac{1}{2}-}}        
\prod_{i=1}^4 \|u_i\|_{X^{0,\frac{1}{2}+}} \\ 
& \lesssim & \frac{N_2^{\frac{1}{4}+} N_4^{\frac{1}{4}+}}{N_1^{\frac{5}{4}-} 
N_3^{\frac{5}{4}-} N} \prod_{i=1}^4 \|u_i\|_{X^{0,\frac{1}{2}+}}
 \lesssim \frac{1 \wedge N_{min}^{0+}}{ N_{max}^{0+} N^{3-}} \prod_{i=1}^4 
\|u_i\|_{X^{0,\frac{1}{2}+}} \, .
\end{eqnarray*}
{\bf b.} $N_1 \ge N \ge N_2$.\\
By the mean value theorem we get
\begin{eqnarray*}
F & \lesssim & \frac{N_2}{N_1} (\frac{N_3}{N})^{\frac{1}{4}}  
(\frac{N_4}{N})^{\frac{1}{4}}   \prod_{i=1}^4 N_i^{-1}  
\frac{N_2^{1+}}{N_1^{\frac{1}{2}-}} \frac{N_4^{1+}}{N_3^{\frac{1}{2}-}}        
\prod_{i=1}^4 \|u_i\|_{X^{0,\frac{1}{2}+}} \\
& \lesssim & \frac{N_2^{1+} N_4^{\frac{1}{4}+}}{N_1^{\frac{5}{2}-} 
N_3^{\frac{5}{4}-} N^{\frac{1}{2}}}  \prod_{i=1}^4 
\|u_i\|_{X^{0,\frac{1}{2}+}}
 \lesssim  \frac{1 \wedge N_{min}^{0+}}{ N_{max}^{0+} N^{3-}} \prod_{i=1}^4 
\|u_i\|_{X^{0,\frac{1}{2}+}} \, .
\end{eqnarray*}
{\bf Case 2:} $N_3 \gtrsim N \ge N_4$ and $N_3 \gg N_4$.\\
{\bf a.} $N_1 \ge N_2 \ge N$.\\
We get similarly as in case 1a.
\begin{eqnarray*}
F & \lesssim & \prod_{i=1}^2(\frac{N_i}{N})^{\frac{1}{4}} \prod_{i=1}^4 
N_i^{-1} 
\frac{N_2^{1+}}{N_1^{\frac{1}{2}-}} \frac{N_4^{1+}}{N_3^{\frac{1}{2}-}}        
\prod_{i=1}^4 \|u_i\|_{X^{0,\frac{1}{2}+}} \\ 
& \lesssim & \frac{N_2^{\frac{1}{4}+} N_4^{0+}}{N_1^{\frac{5}{4}-} 
N_3^{\frac{3}{2}-} N^{\frac{1}{2}}} \prod_{i=1}^4 
\|u_i\|_{X^{0,\frac{1}{2}+}}\\
& \lesssim & \frac{1 \wedge N_{min}^{0+}}{ N_{max}^{0+} N^{3-}} \prod_{i=1}^4 
\|u_i\|_{X^{0,\frac{1}{2}+}} \, .
\end{eqnarray*}
{\bf b.} $N_1 \ge N \ge N_2$.\\
By the mean value theorem we get similarly as in case 1a:
\begin{eqnarray*}
F & \lesssim & \frac{N_2}{N_1} \prod_{i=1}^4 N_i^{-1}  
\frac{N_2^{1+}}{N_1^{\frac{1}{2}-}} \frac{N_4^{1+}}{N_3^{\frac{1}{2}-}}  
\prod_{i=1}^4 \|u_i\|_{X^{0,\frac{1}{2}+}} \\ 
& \lesssim & \frac{N_2^{1+} N_4^{0+}}{N_1^{\frac{5}{2}-} N_3^{\frac{3}{2}-} } 
\prod_{i=1}^4 \|u_i\|_{X^{0,\frac{1}{2}+}}
 \lesssim  \frac{1 \wedge N_{min}^{0+}}{ N_{max}^{0+} N^{3-}} \prod_{i=1}^4 
\|u_i\|_{X^{0,\frac{1}{2}+}} \, .
\end{eqnarray*}
{\bf Case 3:} $N \gg N_3 \ge N_4$ $\Rightarrow N_1 \sim N_2 \gtrsim N$, 
because 
$\sum_{i=1}^4 \xi_i = 0$.\\
We get by (\ref{B2}):
\begin{eqnarray*}
F & \lesssim & (\frac{N_1}{N})^{\frac{1}{4}}(\frac{N_2}{N})^{\frac{1}{4}} 
\prod_{i=1}^4 N_i^{-1} \|u_1 u_3\|_{L_{xt}^2}\|u_2 u_4\|_{L_{xt}^2}  \\
& \lesssim & (\frac{N_1}{N})^{\frac{1}{4}}(\frac{N_2}{N})^{\frac{1}{4}} 
\prod_{i=1}^4 N_i^{-1} \frac{N_3}{N_1^{\frac{1}{2}}} 
\frac{N_4^{1+}}{N_2^{\frac{1}{2}-}}  \prod_{i=1}^4 
\|u_i\|_{X^{0,\frac{1}{2}+}} 
\\ & \lesssim & \frac{ N_4^{0+}}{N_1^{\frac{5}{4}}N_2^{\frac{5}{4}-} 
N^{\frac{1}{2}}} \prod_{i=1}^4 \|u_i\|_{X^{0,\frac{1}{2}+}}
\lesssim  \frac{1 \wedge N_{min}^{0+}}{ N_{max}^{0+} N^{3-}} \prod_{i=1}^4 
\|u_i\|_{X^{0,\frac{1}{2}+}} \, .
\end{eqnarray*}
{\bf 8.} Finally we prove
$$ \int_0^{\delta} |\langle I(u^2) - (Iu)^2, Iu \rangle | dt \lesssim 
N^{-\frac{5}{2}+} \delta^{\frac{1}{2}}\|\nabla Iu\|_{X^{0,\frac{1}{2}+}}^3 \, 
, 
$$
which follows from
\begin{eqnarray}
\label{19}
G  =  |\int_0^{\delta} \int_* M(\xi_1,...,\xi_4) \prod_{i=1}^3 
\widehat{u_i}(\xi_i,t)d\xi_1 d\xi_2 d\xi_3 dt| \lesssim  N^{-\frac{5}{2}+} 
\delta^{\frac{1}{2}} \prod_{i=1}^3 \|u_i\|_{X^{0,\frac{1}{2}+}} \, , 
\end{eqnarray}
where
$$ M(\xi_1,\xi_2,\xi_3) := 
\frac{|m(\xi_1+\xi_2)-m(\xi_1)m(\xi_2)|}{m(\xi_1)m(\xi_2)} 
\prod_{i=1}^3 |\xi_i|^{-1} \, . $$
Assume w.l.o.g. $N_1\ge  N_2$ and $N_1 \ge N$.\\
{\bf Case 1:} $ N_3 \gtrsim N$.\\
{\bf a.} $N_1 \ge N_2 \gtrsim N$.\\
By Strichartz and Sobolev we get:
\begin{eqnarray*}
G & \lesssim & (\frac{N_1}{N})^{\frac{1}{4}}(\frac{N_2}{N})^{\frac{1}{4}} 
\prod_{i=1}^3 N_i^{-1} \|u_1\|_{L_t^2 L_x^6}\| u_2\|_{L_t^{\infty} L_x^3 
}\|u_3\|_{L_{xt}^2}  \\
& \lesssim &  (\frac{N_1}{N})^{\frac{1}{4}}(\frac{N_2}{N})^{\frac{1}{4}} 
\prod_{i=1}^3 N_i^{-1} N_2^{\frac{1}{2}} \delta^{\frac{1}{2}} \prod_{i=1}^3 
\|u_i\|_{X^{0,\frac{1}{2}+}} \\ & \lesssim & \frac{1}{N_1^{\frac{3}{4}}N_3 
N_2^{\frac{1}{4}} N^{\frac{1}{2}}}\delta^{\frac{1}{2}} \prod_{i=1}^3 
\|u_i\|_{X^{0,\frac{1}{2}+}}
 \lesssim  \frac{1 \wedge N_{min}^{0+}}{ N_{max}^{0+} N^{\frac{5}{2}-}} 
\delta^{\frac{1}{2}} \prod_{i=1}^3 \|u_i\|_{X^{0,\frac{1}{2}+}} \, .
\end{eqnarray*}
{\bf b.} $N_1 \ge N \ge N_2$ and $N_1 \gg N_2$.\\
Using the mean value theorem we get
\begin{eqnarray*}
G & \lesssim & \frac{N_2}{N_1} \prod_{i=1}^3 N_i^{-1} N_2^{\frac{1}{2}} 
\delta^{\frac{1}{2}} \prod_{i=1}^3 \|u_i\|_{X^{0,\frac{1}{2}+}} \\
 & \lesssim & \frac{N_2^{\frac{1}{2}}}{N_1^2 N_3} \delta^{\frac{1}{2}} 
\prod_{i=1}^3 \|u_i\|_{X^{0,\frac{1}{2}+}}
 \lesssim  \frac{1 \wedge N_{min}^{0+}}{ N_{max}^{0+} N^{\frac{5}{2}-}}  
\delta^{\frac{1}{2}} \prod_{i=1}^3 \|u_i\|_{X^{0,\frac{1}{2}+}} \, .
\end{eqnarray*}
{\bf Case 2:} $N \gg N_3$ $\Rightarrow N_1 \sim N_2 \gtrsim N$, because 
$\sum_{i=1}^3 = 0$.\\
By (\ref{B2}) we get the estimate
\begin{eqnarray*}
G & \lesssim & (\frac{N_1}{N})^{\frac{1}{4}}(\frac{N_2}{N})^{\frac{1}{4}} 
\prod_{i=1}^3 N_i^{-1} \|u_1 u_3\|_{L_{tx}^3} \|u_2\|_{L_{xt}^2}  \\
& \lesssim &  (\frac{N_1}{N})^{\frac{1}{4}}(\frac{N_2}{N})^{\frac{1}{4}} 
\prod_{i=1}^3 N_i^{-1}\frac{ N_3^{1+}}{N_1^{\frac{1}{2}-}} 
\delta^{\frac{1}{2}} 
\prod_{i=1}^3 \|u_i\|_{X^{0,\frac{1}{2}+}} \\ & \lesssim & 
\frac{N_3^{0+}}{N_1^{\frac{5}{4}-} N_2^{\frac{3}{4}} 
N^{\frac{1}{2}}}\delta^{\frac{1}{2}} \prod_{i=1}^3 
\|u_i\|_{X^{0,\frac{1}{2}+}}
 \lesssim  \frac{1 \wedge N_{min}^{0+}}{ N_{max}^{0+} N^{\frac{5}{2}-}} 
\delta^{\frac{1}{2}} \prod_{i=1}^3 \|u_i\|_{X^{0,\frac{1}{2}+}} \, .
\end{eqnarray*}
This completes the estimates for the increment of $E(Iu)$ on the local 
existence 
interval $[0,\delta]$ in terms of the parameter $N$. 

We recall our aim to give an a-priori bound of $\|\nabla Iu(t)\|_{L^2}$ (cf. 
(\ref{6})) on any interval $[0,T]$. We want to show this as a consequence of 
Proposition \ref{LWP1} and the estimates for the modified energy just given.

We assume $N\ge 1$ to be a number to be specified later and $s\ge 
\frac{3}{4}$. 
Let data $u_0 \in H^s({\mathbb R}^3)$ be given. Then we have
\begin{eqnarray}
  \label{21}
\|\nabla Iu_0\|_{L^2}^2 & \lesssim& \| |\xi|\widehat{u_0}(\xi)\|_{L^2(\{|\xi| 
\le 
N\})}^2 + \|N^{1-s} |\xi|^s \widehat{u_0}(\xi)\|_{L^2(\{|\xi| \ge N\})}^2  \\ 
\nonumber 
& \lesssim &  \|N^{1-s} |\xi|^s \widehat{u_0}(\xi)\|_{L^2({\mathbb R}^3)}^2 = 
N^{2(1-s)} \|u_0\|_{\dot{H}^s}^2 \lesssim N^{2(1-s)} \, . 
\end{eqnarray}
This immediately implies an estimate for $E(Iu_0)$. We namely have for $s\ge 
\frac{3}{4}$:
$$ \|Iu_0\|_{L^4({\mathbb R}^3)}^4 \lesssim \|Iu_0\|_{\dot{H}^{\frac{3}{4}}}^4 
\lesssim \|u_0\|_{H^s}^4 $$ and trivially $\|Iu_0\|_{L^2}^2 \lesssim 
\|u_0\|_{L^2}^2$. Thus using the definition of the modified energy and 
(\ref{21}):
$$ E(Iu_0) \le c_0 N^{2(1-s)} \, . $$
From this we get
$$ \|\nabla Iu_0\|_{L^2}^2 \le c_0 N^{2(1-s)} \, . $$
Our local existence theorem (Proposition \ref{LWP1}) shows that the Cauchy 
problem 
(\ref{0.3}),(\ref{0.4}) has a unique solution $u$ with $\nabla Iu \in 
X^{0,\frac{1}{2}+}[0,\delta]$ and $$  \|\nabla Iu(\delta)\|_{L^2}^2 \le 
\|\nabla 
Iu\|_{X^{0,\frac{1}{2}+}[0,\delta]}^2 \le 2 \|\nabla Iu_0\|_{L^2}^2 \le 2c_0 
N^{2(1-s)} \, . $$
Here $\delta$ can be chosen subject to the conditions (\ref{Delta}) , namely 
(because $s > 1/2$): 
$$ \max(\delta^{s-\frac{1}{2}-},\delta^{\frac{s}{2}-} 
N^{1-s},\delta^{\frac{1}{2}-} N^{2(1-s)}) \sim 1 \Longleftrightarrow \delta 
\sim 
  \frac{1}{N^{4(1-s)+}} \, .$$
In order to reapply the local existence theorem with time steps of equal 
length 
we need a uniform bound of $\|\nabla Iu(t)\|_{L^2}$ at time $t=\delta$, 
$t=2\delta$ 
etc., which follows from a uniform control over the modified energy. The 
increment of the energy is controlled by (\ref{E}) and the estimates of this 
section as follows, provided $s \ge 3/4$:
\begin{eqnarray*}
\lefteqn{ |E(Iu(\delta))-E(Iu_0)|}\\
&\hspace{-1em} \lesssim &N^{-1+} \|\nabla 
Iu\|_{X^{0,\frac{1}{2}+}(0,\delta)}^4 
+ N^{-1+} \delta^{\frac{1}{4}-}\|\nabla Iu\|_{X^{0,\frac{1}{2}+}(0,\delta)}^3 
+ 
N^{-2+} \|\nabla Iu\|_{X^{0,\frac{1}{2}+}(0,\delta)}^6 \\
& & + N^{-\frac{5}{2}+} \|\nabla Iu\|_{X^{0,\frac{1}{2}+}(0,\delta)}^5 + 
N^{-3+} 
\|\nabla Iu\|_{X^{0,\frac{1}{2}+}(0,\delta)}^4 + N^{-\frac{5}{2}+} 
\delta^{\frac{1}{2}} \|\nabla Iu\|_{X^{0,\frac{1}{2}+}(0,\delta)}^3  
\end{eqnarray*}
The last two terms can be neglected in comparison to the others. Thus we get
\begin{eqnarray*}
\lefteqn{ |E(Iu(\delta))-E(Iu_0)|}\\
& \lesssim &N^{-1+} N^{4(1-s)} + N^{-1+}N^{-(1-s)+} N^{3(1-s)} + N^{-2+} 
N^{6(1-s)} + N^{-\frac{5}{2}+} N^{5(1-s)} \, .  
\end{eqnarray*}
One easily checks that the first term is the decisive one, so that
$$ |E(Iu(\delta))-E(Iu_0)| \le c_1 N^{-1+} N^{4(1-s)} \, ,$$
where $c_1 = c_1(c_0)$. This is easily seen to be bounded by $c_0 N^{2(1-s)}$ 
for large $N$. The number of iteration steps to reach a given time $T$ is 
$\frac{T}{\delta} \sim TN^{4(1-s)+}$. This means that in order to give a 
uniform 
bound of the energy of the iterated solutions, namely by $2c_0 N^{2(1-s)}$, 
from 
the last inequality, the following condition has to be fulfilled:
$$ c_1 N^{-1+} N^{4(1-s)} T N^{4(1-s)+} < c_0 N^{2(1-s)} \, , $$
where $c_1 = c_1(2c_0)$ (recall here that the initial energy is bounded by 
$c_0
N^{2(1-s)}$). This can be fulfilled for $N$ sufficiently large, provided
$$ -1 +4(1-s) + 4(1-s) < 2(1-s) \Longleftrightarrow s > 5/6 \, . $$
This gives the desired bound for $\|\nabla Iu(t)\|_{L^2}$ on any interval 
$[0,T]$, i.e. (\ref{6}) is proved. As explained above this completes the proof 
of our Theorem \ref{Theorem 1}.


\begin{thebibliography}{99999999}
\bibitem[BS]{BS} F. Bethuel and J.C. Saut: {\sl Travelling waves for the 
Gross-Pitaevskii equation I}. Ann. I. H. Poincar\'e Phys. Th\'eor. 70 (1999), 
147-238
\bibitem[B]{B} J. Bourgain: {\sl Scattering in the energy space and below for 
3D 
NLS}. J. d'Analyse Math. 75 (1998), 267-297
\bibitem[CH]{CH} T. Cazenave and A. Haraux: {\sl An introduction to semilinear 
evolution equations}. Oxford science publications 1998 
\bibitem[CW]{CW} T. Cazenave and F. Weissler: {\sl The Cauchy problem for the 
nonlinear Schr\"odinger equation in $H^s$}. Nonlinear Analysis 14 (1990), 
807-836
\bibitem[CKSTT]{CKSTT} J. Colliander, M. Keel, G. Staffilani, H. Takaoka and 
T. Tao: 
{\sl Almost conservation laws and global rough solutions to a nonlinear 
Schr\"odinger equation}. Math. Res. Letters 9 (2002), 659-682 
\bibitem[Ga]{Ga} C. Gallo: {\sl The Cauchy problem for defocusing nonlinear 
Schr\"odinger equations with non-vanishing initial data at infinity}. Comm. 
Part. Diff. Equa. 33 (2008), 729-771
\bibitem[Ge]{Ge} P. G\'erard: {\sl The Cauchy problem for the Gross-Pitaevskii 
equation}. Ann. I. H. Poincar\'e Anal. Non-lin\'eaire 23 (2006), 765-779
\bibitem[GTV]{GTV} J. Ginibre, Y. Tsutsumi and G. Velo: {\sl On the Cauchy 
problem for the Zakharov system}. J. Funct. Analysis 151 (1997), 384-436
\bibitem[Gr]{Gr} E.P. Gross: {\sl Hydrodynamics of a Superfluid Condensate}. 
J. 
Math. Phys. 4 (1963), 195-207
\bibitem[G]{G} A. Gr\"unrock: {\sl New applications of the Fourier restriction 
norm method to wellposedness problems for nonlinear evolution equations}. 
Dissertation Univ. Wuppertal 2002, \\ 
{\tt http://elpub.bib.uni-wuppertal.de/servlets/DocumentServlet?id=254}
\bibitem[K]{K} T. Kato: {\sl On nonlinear Schr\"odinger equations II. 
$H^s$-solutions and unconditional well-posedness}. J. d'Analyse Math. 67 
(1995), 
281-306
\bibitem[KT]{KT} M. Keel and T. Tao: {\sl Endpoint Strichartz estimates}. 
Amer. 
J. Math. 120 (1998), 955-980
\bibitem[P]{P} L.P. Pitaevskii: {\sl Vortex lines in an imperfect Bose gas}. 
Soviet Physics JETP 13 (1961), 451-454
\bibitem[SS]{SS} C. Sulem and P.L. Sulem: {\sl The nonlinear Schr\"odinger 
equation. Self-focusing and wave collapse}. Appl. Math. Sci. vol. 139, 
Springer 
1999
\end{thebibliography}
\end{document}